\documentclass[11pt,a4paper]{article}
\usepackage{amsfonts}
\usepackage{latexsym}
\usepackage{cite}
\usepackage{amsmath,amsfonts,latexsym,amssymb}
\usepackage[mathscr]{eucal}
\usepackage{cases,color}
\usepackage{amsthm}
\usepackage{hyperref}
\usepackage{makecell}
\usepackage{multirow}
\usepackage[bf,small]{caption2}
\usepackage{float}
\usepackage{graphicx}
\usepackage{amsmath}
\usepackage{amssymb}
\usepackage[all]{xy}

\newtheorem{theorem}{theorem}[section]
\newtheorem{thm}[theorem]{Theorem}
\newtheorem{lem}[theorem]{Lemma}

\newtheorem{prop}[theorem]{Proposition}
\newtheorem{cor}[theorem]{Corollary}

\newtheorem{defn}[theorem]{Definition}
\newtheorem{exmp}[theorem]{Example}

\newtheorem{rmk}[theorem]{Remark}

\newtheorem{conv}[theorem]{Convention}
\DeclareMathOperator\im{im}

\begin{document}

\title{\vspace{-2cm}\textbf{Finite $p$-groups of class $2$ as central extensions}}
\author{\Large Haimiao Chen}

\date{}
\maketitle

\begin{abstract}
  Finite $p$-groups of nilpotency class 2 are treated from the perspective of central extensions.
  Given finite abelian groups $G,A$, we derive an explicit formula for cocycles representing elements of $H^2(G,A)$, compute $H^2(G,A)$, and describe the actions of ${\rm End}(G)$ and ${\rm End}(A)$ on $H^2(G,A)$. These are used to provide an efficient criterion for lifting endomorphisms of $G$ to homomorphisms between two central extensions. Subsequently, we present two applications to illustrate the usefulness of this approach, in the case $p>2$.
  First, we recover the classification of two-generator $p$-groups of class $2$ up to isomorphism, and compute the order of the automorphism group for each isomorphism class. Second, we construct a family of nonabelian $p$-groups of order $p^7$ whose automorphism groups are abelian.

\end{abstract}

\section{Introduction}

Throughout, $p$ always denotes a prime number.
Finite $p$-groups form a rich variety of objects and attracted much attention of literature. There are many interesting problems concerned with various topics.

In particular, automorphisms and automorphism groups of $p$-groups have been widely studied (see \cite{Ab07,AAG17,FO14,GJ15,HM06,JRY13,JY12,KY18,LW10,Ya13,Wi72} and the references therein).
Until now, automorphism groups have been investigated in-depth for only a few families of groups, including abelian groups \cite{HR07}, metacyclic groups \cite{BC06,Cu07,Cu08,CXZ18}, direct products \cite{Bi08,Br25,BCM06}, certain central extensions of $\mathbb{Z}_{p^k}^n$ by $\mathbb{Z}_{p^k}$\cite{WL21}, a class of $3$-generator $p$-groups \cite{Sz22}, and so on.

Among nonabelian finite $p$-groups, those of nilpotency class $2$ can be considered as the most accessible ones. However, they are still far from well-understood.

In recent decades, an interesting specific topic was formed by {\it Miller groups} \cite{Ca15,JRY13,JY12,KY18,Mo94}, which by definition are nonabelian groups with abelian automorphism groups.
Much state-of-the-art information is contained in the survey \cite{KY18}. It is known that Miller groups are necessarily of class $2$, and a Miller $p$-group has order at least $p^7$ (resp. $2^6$) if $p>2$ (resp. $p=2$).

To study $p$-groups of class $2$, we take a very natural approach.

If $L$ is a $p$-group of class 2, then $L$ is a central extension of the abelianization $L^{\rm ab}$ by the commutator subgroup $[L,L]$, and
$L^{\rm ab}$, $[L,L]$ are nontrivial.

On the other hand, given nontrivial abelian $p$-groups $A,G$, and an $A$-valued 2-cocycle $\alpha$ on $G$, i.e., a map $\alpha:G\times G\to A$ satisfying
\begin{align*}
\alpha(z_1,z_2)-\alpha(z_1,z_2z_3)+\alpha(z_1z_2,z_3)-\alpha(z_2,z_3)=0
\end{align*}
for all $z_1,z_2,z_3\in G$, one can construct a group $E(\alpha)$ with underlying set $A\times G$ and group operation defined by
\begin{align}
(a_1,z_1)\cdot(a_2,z_2)=(a_1+a_2+\alpha(z_1,z_2),z_1z_2).  \label{eq:product}
\end{align}
Here we denote the group operation in $A$ additively, and denote that in $G$ multiplicatively, as done throughout.

As a basic fact in homological algebra (see \cite[Theorem 10.3]{HS97} for instance), up to certain equivalence the central extension $1\to A\to E(\alpha)\to G\to 1$ is uniquely determined by the class $[\alpha]\in H^2(G,A)$, where the $\mathbb{Z}G$-module structure on $A$ is trivial. Furthermore, each class has a {\it normalized} representative $\alpha$, meaning that $\alpha(1,z)=\alpha(z,1)=0$ for all $z\in G$.

From (\ref{eq:product}) we see that $(\alpha(z_1,z_2),z_1z_2)=(0,z_1)\cdot(0,z_2)$ for any $z_1,z_2$. So $\alpha$ can be recovered from $E(\alpha)$.

According to (\ref{eq:product}), in $E(\alpha)$ we have
\begin{align}
(a,z)^{-1}&=(-a-\alpha(z,z^{-1}),z^{-1}), \nonumber \\
[(a_1,z_1),(a_2,z_2)]&=(\eta^\alpha(z_1,z_2),1),  \label{eq:commutator}
\end{align}
where the commutator $[u,v]=uvu^{-1}v^{-1}$, and
\begin{align}
\eta^\alpha(z_1,z_2)=\alpha(z_1,z_2)-\alpha(z_2,z_1).  \label{eq:eta-def}
\end{align}
By (\ref{eq:commutator}), $[E(\alpha),E(\alpha)]=A$ if and only if $\im(\eta^\alpha)$ generates $A$.

Let $Z^2_\circ(G,A)$ denote the set of normalized 2-cocycles $\alpha$ such that $\im(\eta^\alpha)$ generates $A$.

The following is basic and standard. We prove it for self-containedness.
\begin{prop} \label{prop:basic}
Suppose $\alpha,\alpha'\in Z^2_\circ(G,A)$.
\begin{enumerate}
  \item[\rm(a)] If $\langle z_1,\ldots,z_n\rangle=G$, then $\langle(a_1,z_1),\ldots,(a_n,z_n)\rangle=E(\alpha)$ for any $a_i\in A$.
  \item[\rm(b)] Each homomorphism $E(\alpha)\to E(\alpha')$ induces an endomorphism of $A$ and an endomorphism of $G$.
  \item[\rm(c)] An endomorphism $\varphi\in{\rm End}(G)$ can be lifted to $E(\alpha)\to E(\alpha')$
         if and only if $\varphi^\ast[\alpha']=\psi_\ast[\alpha]$ in $H^2(G,A)$, for some $\psi\in{\rm End}(A)$.

         When this holds, each lifting is given by $(a,z)\mapsto(\psi(a)+\mu(z),\varphi(z))$ for a unique map $\mu:G\to A$ with $\partial\mu=\varphi^\ast\alpha'-\psi_\ast\alpha$, where $\partial\mu$ sends $(z_1,z_2)$ to $\mu(z_1z_2)-\mu(z_1)-\mu(z_2)$.
\end{enumerate}
\end{prop}

\begin{proof}
(a) For each $z\in G$, take $k_1,\ldots,k_n\in\mathbb{Z}$ such that $z=z_1^{k_1}\cdots z_n^{k_n}$, and write $(a_1,z_1)^{k_1}\cdots(a_n,z_n)^{k_n}$ as $(b_z,z)$, with $b_z\in A$. For any $z,z'$, by (\ref{eq:commutator}),
$$(\eta^\alpha(z,z'),1)=[(b_z,z),(b_{z'},z')]\in\langle(a_1,z_1),\ldots,(a_n,z_n)\rangle=:K.$$
Hence $A=\langle\im(\eta^\alpha)\rangle$ is identified with a subgroup of $K$ via $a\mapsto(a,1)$. Since $(a,z)=(a-b_z,1)\cdot(b_z,z)\in K$ for all $a\in A$, $z\in G$, we have $E(\alpha)=K$.

\medskip

(b) This is simply due to the fact that $[E(\alpha),E(\alpha)]=[E(\alpha'),E(\alpha')]=A$ and $E(\alpha)^{\rm ab}=E(\alpha')^{\rm ab}=G$,

\medskip

(c) Suppose $\phi$ is a lifting of $\varphi$. Let $\psi=\phi|_A$. Then there exists a unique map $\mu:G\to A$ such that
$\phi(0,z)=(\mu(z),\varphi(z))$ for all $z\in G$.
Since $\alpha$ is normalized, we have $(a,1)(0,z)=(a,z)$, so that
\begin{align*}
\phi(a,z)=\phi(a,1)\phi(0,z)=(\psi(a),1)\cdot(\mu(z),\varphi(z))=(\psi(a)+\mu(z),\varphi(z)),
\end{align*}
the last equality owing to the fact that $\alpha'$ is normalized.

The condition $\phi$ being a homomorphism is equivalent to
\begin{align*}
\alpha'(\varphi(z_1),\varphi(z_2))-\psi(\alpha(z_1,z_2))=\mu(z_1z_2)-\mu(z_1)-\mu(z_2),
\end{align*}
i.e., $\varphi^\ast\alpha'-\psi_\ast\alpha=\partial\mu$. So $\varphi^\ast[\alpha']=\psi_\ast[\alpha]$ in $H^2(G,A)$.

Conversely, suppose $\varphi^\ast[\alpha']=\psi_\ast[\alpha]$ for some $\psi\in{\rm End}(A)$, then there exists $\mu:G\to A$ with $\varphi^\ast\alpha'-\psi_\ast\alpha=\partial\mu$. From the previous paragraph we see that the map
$(a,z)\mapsto(\psi(a)+\mu(z),\varphi(z))$ is a lifting of $\varphi$.

The proof is completed.
\end{proof}

\begin{rmk}\label{rmk:determine}
\rm In the context of Proposition \ref{prop:basic} (c), actually $\psi$ is determined by $\varphi$, since $\langle\im(\eta^\alpha)\rangle=A$ and for any $z_1,z_2$,
\begin{align*}
(\psi(\eta^{\alpha}(z_1,z_2)),1)
&=\phi([(0,z_1),(0,z_2)])=[(\mu(z_1),\varphi(z_1)),(\mu(z_2),\varphi(z_2))]  \\
&=\big(\eta^{\alpha'}(\varphi(z_1),\varphi(z_2)),1\big).
\end{align*}
\end{rmk}

\medskip

In Section 2, we derive an explicit formula for 2-cocycles representing elements of $H^2(G,A)$, and work out the actions of ${\rm End}(G)$ and ${\rm End}(A)$ on $H^2(G,A)$. They have independent interests. Based on these, in Section 3, assuming $p>2$, we present an efficient criterion for lifting endomorphisms of $G$ to homomorphisms between two central extensions. We expect the results of Section 2 and Section 3 to be potentially applicable in dealing with $p$-groups of class $2$ as well as their automorphism groups (or endomorphism semigroups). Automorphism groups of group extensions fit in the famous {\it Wells exact sequence} \cite{We71}, and were understood well \cite{JL10,PSY10}. The explicit knowledge on the second cohomology group will make the Wells sequence more useful. In Section 4, as an application, for $p>2$, we recover the result of \cite{AMM12} which classifies two-generator $p$-groups of class $2$ up to isomorphism, and determine the order of the automorphism group for each isomorphism class.
In Section 5, as another application, for $p>2$, we construct a family of Miller groups of order $p^7$. Previously, only one such group was known.

\medskip

\begin{conv}
\rm Use $\mathbb{Z}_{>0}$ to denote the set of positive integers.

For a ring $R$, write elements of $R^k$ as column vectors.

Given $\ell,h\in\mathbb{Z}_{>0}$ with $h>\ell$, each $x\in\mathbb{Z}$ is regarded as an element of $\mathbb{Z}_{p^\ell}$ via the canonical quotient map $\mathbb{Z}\twoheadrightarrow \mathbb{Z}_{p^\ell}$, and each $x\in\mathbb{Z}_{p^{h}}$ is regarded as an element of $\mathbb{Z}_{p^\ell}$ via the reduction map $\mathbb{Z}_{p^{h}}\twoheadrightarrow \mathbb{Z}_{p^\ell}$.

Let $\Lambda=\{(i,j)\colon 1\le i<j\le n\}$. Once an enumeration of elements of $\Lambda$ is chosen, $\Lambda$ will be identified with $\{1,\ldots,{n\choose 2}\}$, through which the $k$-th element of $\Lambda$ corresponds to $k$.
Then for any set $X$, the Cartesian product $X^{\Lambda}$ will be identified with $X^{{n\choose 2}}$.
\end{conv}

Here are some notations which will be adopted throughout:
\begin{center}
\begin{tabular}{|c|c|c|}
  \hline
  Symbol & Meaning & Setting/Instruction  \\
  \hline
  $\mathcal{M}_k^\ell(R)$ &\makecell{the set of $k\times \ell$ matrices with \\ entries in $R$} & $k,\ell\in\mathbb{Z}_{>0}$, $R$ is a ring \\
  \hline
  $\mathcal{M}_k(R)$ & a shorthand for $\mathcal{M}_k^k(R)$ &  \\
  \hline
  $\mathbf{x}^{\rm tr}$ & the transpose of $\mathbf{x}$ & \multirow{3}{*}{$\mathbf{x}=(x_{ij})$ is a matrix}   \\
  $\mathbf{x}_{i}$  &  the $i$-th row of $\mathbf{x}$ &  \\
  $\mathbf{x}_{i,j}$  &$x_{ij}$  &  \\
  \hline
  $|X|$ & the cardinality of $X$ & $X$ is a finite set   \\
  \hline
  $\|c\|$ & the largest integer $e\le\ell$ with $p^e\mid c$ & $c\in \mathbb{Z}_{p^\ell}$, with $\ell\in\mathbb{Z}_{>0}$  \\
  \hline
  $G$ & $\mathbb{Z}_{p^{m_1}}\times\cdots\times \mathbb{Z}_{p^{m_n}}$ & \makecell{the $m_i$'s are fixed, with \\ $1\le m_1\le\cdots\le m_n:=m$}  \\
  \hline
  $g_i$ & a fixed generator of the factor $\mathbb{Z}_{p^{m_i}}$ &   \\
  \hline
  $\mathbf{g}^{\mathsf{u}}$ & the element $g_1^{u_1}\cdots g_n^{u_n}\in G$ &  $\mathsf{u}=(u_1,\ldots,u_n)^{\rm tr}\in\mathbb{Z}^n$  \\
  \hline
  $A$ & $\mathbb{Z}_{p^{d_1}}\times\cdots\times \mathbb{Z}_{p^{d_r}}$ & \makecell{the $d_i$'s are fixed, with \\ $1\le d_1\le\cdots\le d_r:=d$}  \\
  \hline
  $\hat{k}$ & \makecell{the endomorphism of $A$ \\ given by $a\mapsto p^ka$}  & $k\in\mathbb{Z}_{>0}$ \\
  \hline
  $\delta_{a,b}$ & \makecell{the Kronecker symbol: $\delta_{a,b}=1$ if \\ $a=b$, and $\delta_{a,b}=0$ if $a\ne b$} & $a,b$ are elements of a set   \\
  \hline
  $[x]_{p^h}$ & \makecell{the unique $z\in\{0,1,\ldots,p^h-1\}$ \\ with $z\equiv x\pmod{p^h}$}
              & \makecell{$x\in\mathbb{Z}$ or $x\in \mathbb{Z}_{p^\ell}$ \\  for some $\ell\ge h$}  \\
  \hline
  $\zeta_{p^h}(x,y)$ & \makecell{$\zeta_{p^h}(x,y)=1$ if $[x]_{p^h}+[y]_{p^h}\ge p^h$,   \\ $\zeta_{p^h}(x,y)=0$ otherwise}
                 & \makecell{$x,y\in\mathbb{Z}$ or $x,y\in \mathbb{Z}_{p^\ell}$ \\  for some $\ell\ge h$}  \\
  \hline
\end{tabular}
\end{center}

\section{$H^2(G,A)$}

\subsection{An explicit formula for $2$-cocycles}

Given $\mathbf{a}=(a_1,\ldots,a_n)\in A^n$ and $\mathbf{c}=\{c_{ij}\}_{(i,j)\in\Lambda}\in A^{\Lambda}$,
define
$$\alpha_{\mathbf{a},\mathbf{c}}:G\times G\to A$$
as follows: for any $\mathsf{s}=(s_1,\ldots,s_n)^{\rm tr}, \mathsf{t}=(t_1,\ldots,t_n)^{\rm tr}\in\mathbb{Z}^n$, put
\begin{align}
\alpha_{\mathbf{a},\mathbf{c}}(\mathbf{g}^{\mathsf{s}},\mathbf{g}^{\mathsf{t}})
=\sum\limits_{i=1}^n\zeta_{p^{m_i}}(s_i,t_i)a_i-\sum\limits_{(i,j)\in\Lambda}t_is_jc_{ij}.   \label{eq:alpha}
\end{align}
It is not difficult to see that
$$\alpha_{\mathbf{a},\mathbf{c}}(\mathbf{g}^{\mathsf{u}},\mathbf{g}^{\mathsf{v}})
-\alpha_{\mathbf{a},\mathbf{c}}(\mathbf{g}^{\mathsf{u}},\mathbf{g}^{\mathsf{v}+\mathsf{w}})
+\alpha_{\mathbf{a},\mathbf{c}}(\mathbf{g}^{\mathsf{u}+\mathsf{v}},\mathbf{g}^{\mathsf{w}})
-\alpha_{\mathbf{a},\mathbf{c}}(\mathbf{g}^{\mathsf{v}},\mathbf{g}^{\mathsf{w}})=0$$
for any $\mathsf{u},\mathsf{v},\mathsf{w}\in\mathbb{Z}^n$, and $\alpha_{\mathbf{a},\mathbf{c}}(\mathbf{g}^{\mathsf{s}},\mathbf{g}^{\mathsf{t}})=0$ when $\mathsf{s}=\mathsf{0}$ or $\mathsf{t}=\mathsf{0}$. Hence $\alpha_{\mathbf{a},\mathbf{c}}$ is a normalized $2$-cocycle.

One can recover $\mathbf{a},\mathbf{c}$ from $\alpha_{\mathbf{a},\mathbf{c}}$ via
$$a_i=\alpha_{\mathbf{a},\mathbf{c}}(g_i,g_i^{-1}), \qquad    c_{ij}=-\alpha_{\mathbf{a},\mathbf{c}}(g_j,g_i).$$

\begin{thm} \label{thm:H2}
Each element of $H^2(G,A)$ is represented by $\alpha_{\mathbf{a},\mathbf{c}}$ for some $\mathbf{a}\in A^n$, $\mathbf{c}\in A^{\Lambda}$ with $c_{ij}\in\ker(\widehat{m_i})$ for all $(i,j)\in\Lambda$.
Moreover, $[\alpha_{\mathbf{a},\mathbf{c}}]=0$ if and only if
$a_i\in\im(\widehat{m_i})$ for all $i\in\{1,\ldots,n\}$ and $c_{ij}=0$ for all $(i,j)\in\Lambda$.

Consequently, the homomorphism
$A^n\times\prod_{(i,j)\in\Lambda}\ker(\widehat{m_i})\to H^2(G,A)$, $(\mathbf{a},\mathbf{c})\mapsto[\alpha_{\mathbf{a},\mathbf{c}}]$
induces the inverse of an isomorphism
\begin{align*}
H^2(G,A)\cong\prod\limits_{i=1}^n{\rm coker}(\widehat{m_i})\times\prod_{(i,j)\in\Lambda}\ker(\widehat{m_i})
\end{align*}
which is induced by
$\alpha\mapsto\big((\alpha(g_1,g_1^{-1}),\ldots,\alpha(g_n,g_n^{-1})),\{-\alpha(g_j,g_i)\}_{(i,j)\in\Lambda}\big).$
\end{thm}

\begin{proof}
We will slightly modify the proof of \cite[Lemma 3.3]{HLY14} and generalize it. The main steps are: (i) compute $2$-cocycles $\beta$ of a Koszul-like resolution $(K_\bullet,d_\bullet)$, and (ii) use a chain map from the normalized bar resolution $(B_\bullet,\partial_\bullet)$ to $(K_\bullet,d_\bullet)$, to pull back each $\beta$ to be a map $G\times G\to A$.

It is well-known that the periodic sequence
$$\cdots\to\mathbb{Z}\mathbb{Z}_{p^{m_i}}\stackrel{T_i}\longrightarrow\mathbb{Z}\mathbb{Z}_{p^{m_i}}\stackrel{N_i}\longrightarrow
\mathbb{Z}\mathbb{Z}_{p^{m_i}}\stackrel{T_i}\longrightarrow\mathbb{Z}\mathbb{Z}_{p^{m_i}}\stackrel{N_i}\longrightarrow\mathbb{Z}\longrightarrow 0$$
is a free resolution of the trivial $\mathbb{Z}\mathbb{Z}_{p^{m_i}}$-module $\mathbb{Z}$,
where $T_i$ is the multiplication by $g_i-1$, and $N_i$ is the multiplication by $\sum_{j=0}^{p^{m_i}-1}g_i^j$.

A free resolution of the trivial $\mathbb{Z}G$-module $\mathbb{Z}$ is given by the tensor product $(K_\bullet,d_\bullet)$ of the above resolutions for the $\mathbb{Z}_{p^{m_i}}$'s.
To be explicit,
$$K_q=\bigoplus\limits_{u_1+\cdots+u_n=q}(\mathbb{Z}G)\Phi_{(u_1,\ldots,u_n)},$$
where $\Phi_{(u_1,\ldots,u_n)}$ denotes a free generator indexed by $(u_1,\ldots,u_n)$ with $u_i\ge 0$,
and the differential is $d=d_1+\cdots+d_n$, with
$$d_i\Phi_{(u_1,\ldots,u_n)}=
\begin{cases}
0, & u_i=0 \\
(-1)^{\sum_{s<i}u_s}N_i\Phi_{(u_1,\ldots,u_i-1,\ldots,u_n)}, & 2\mid u_i\ne0  \\
(-1)^{\sum_{s<i}u_s}T_i\Phi_{(u_1,\ldots,u_i-1,\ldots,u_n)}, & 2\nmid u_i\ne0 \\
\end{cases}.$$

Let $\mathsf{e}^i=(0,\ldots,1,\ldots,0)$, the $n$-component vector in which the unique $1$ sits at the $i$-th position.
Introduce
\begin{alignat*}{3}
\rho_i&=\Phi_{\mathsf{e}^i},  \qquad &\sigma_i&=\Phi_{2\mathsf{e}^i}, \qquad  &\sigma_{ij}&=\Phi_{\mathsf{e}^i+\mathsf{e}^j}, \\
\tau_i&=\Phi_{3\mathsf{e}^i}, \qquad &\tau_{ij}&=\Phi_{2\mathsf{e}^i+\mathsf{e}^j},\qquad &\tau_{ijk}&=\Phi_{\mathsf{e}^i+\mathsf{e}^j+\mathsf{e}^k}.
\end{alignat*}
A 2-cochain $\beta\in\hom_{\mathbb{Z}G}(K_2,A)$ is determined by $\beta(\sigma_i)$, $1\le i\le n$, and $\beta(\sigma_{ij})$, $1\le i<j\le n$.
We shall remember that $A$ is a trivial $\mathbb{Z}G$ module.

By definition, $\beta$ is a cocycle if and only if
$$0=(d^\ast\beta)(\tau_{i_1\cdots i_s})=\beta(d\tau_{i_1\cdots i_s})$$
for all $1\le i_1<\cdots<i_s\le n$ with $1\le s\le 3$.
Note that
\begin{align*}
\beta(d\tau_i)&=\beta(T_i\sigma_i)=0,    \\
\beta(d\tau_{ij})&=\beta(T_j\sigma_i\pm N_i\sigma_{ij})=\pm p^{m_i}\beta(\sigma_{ij}),  \\
\beta(d\tau_{ijk})&=\beta(T_i\sigma_{jk}-T_j\sigma_{ik}+T_k\sigma_{ij})=0.
\end{align*}
Hence $\beta$ is a cocycle if and only if $p^{m_i}\beta(\sigma_{ij})=0$ for all $i<j$.

If $\beta=d^\ast\gamma$ for some 1-cochain $\gamma$, then
\begin{alignat*}{2}
\beta(\sigma_i)&=\gamma(d\sigma_i)=p^{m_i}\gamma(\rho_i) \qquad &&\text{for\ all\ }i, \\
\beta(\sigma_{ij})&=\gamma(d\sigma_{ij})=0  \qquad &&\text{for\ all\ } i<j.
\end{alignat*}
Conversely, if $\beta(\sigma_i)=p^{m_i}c_i$ with $c_i\in A$ for all $i$, and $\beta(\sigma_{ij})=0$ for all $(i,j)\in\Lambda$, then $\beta=d^\ast\gamma$, where $\gamma$ is determined by $\rho_i\mapsto c_i$.

Plug in the map $F_2:B_2\to K_2$ given on \cite[Page 4]{HWY20}, namely,
\begin{align*}
[\mathbf{g}^{\mathsf{s}},\mathbf{g}^{\mathsf{t}}]\mapsto&\sum_{i=1}^ng_1^{s_1+t_1}\cdots g_{i-1}^{s_{i-1}+t_{i-1}}\Big[\frac{s_i+t_i}{p^{m_i}}\Big]
\sigma_i   \\
&-\sum_{1\le i<j\le n}g_1^{s_1}\cdots g_{j-1}^{s_{j-1}}g_1^{t_1}\cdots g_{i-1}^{t_{i-1}}\left({\sum}_{k=0}^{t_i-1}g_i^k\right)\left({\sum}_{k=0}^{s_j-1}g_j^k\right)\sigma_{ij},
\end{align*}
for $0\le s_i,t_i<p^{m_i}$.
Let $a_i=\beta(\sigma_i)$, $c_{ij}=\beta(\sigma_{ij})$, then
$$(F_2^\ast\beta)([\mathbf{g}^{\mathsf{s}},\mathbf{g}^{\mathsf{t}}])
=\sum\limits_{i=1}^n\zeta_{p^{m_i}}(s_i,t_i)a_i-\sum\limits_{1\le i<j\le n}t_is_jc_{ij}.$$

The proof is completed.
\end{proof}

\begin{exmp}
\rm Suppose $r=2$, $n=3$, $d_1<m_1<d_2<m_2<m_3.$
Then
\begin{align*}
{\rm coker}(\widehat{m_1})\cong\ker(\widehat{m_1})\cong \mathbb{Z}_{p^{d_1}}\times \mathbb{Z}_{p^{m_1}}=:A_0, \\
{\rm coker}(\widehat{m_i})=\ker(\widehat{m_i})=A, \qquad i=2,3.
\end{align*}
Hence
$$H^2(G,A)\cong (A_0\times A^2)\times (A_0^2\times A).$$
\end{exmp}

\begin{rmk}\label{rmk:invariants}
\rm From the definition (\ref{eq:eta-def}) we see
\begin{align}
\eta^{\alpha_{\mathbf{a},\mathbf{c}}}(\mathbf{g}^{\mathsf{s}},\mathbf{g}^{\mathsf{t}})
=\sum_{1\le i<j\le n}(s_it_j-s_jt_i)c_{ij}, \label{eq:eta}
\end{align}
so the condition $\alpha_{\mathbf{a},\mathbf{c}}\in Z^2_\circ(G,A)$, i.e. $\langle\im(\eta^{\alpha_{\mathbf{a},\mathbf{c}}})\rangle=A$,
is equivalent to
\begin{align}
\langle c_{12},\ldots,c_{n-1,n}\rangle=A.   \label{eq:generate}
\end{align}

A straightforward consequence is $r\le{n\choose 2}$.
Moreover, since $c_{ij}\in\ker(\widehat{m_{i}})\leqslant\ker(\widehat{m_{n-1}})$ for all $(i,j)\in\Lambda$, we have $p^{m_{n-1}}A=0$, so
$m_{n-1}\ge d$.
\end{rmk}

\begin{defn}
\rm Call a pair $(\mathbf{a},\mathbf{c})\in A^n\times A^\Lambda$ {\it admissible} if $\mathbf{c}=\{c_{ij}\}_{(i,j)\in\Lambda}$ satisfies $c_{ij}\in\ker(\widehat{m_i})$ for all $(i,j)$ and $\langle c_{12},\ldots,c_{n-1,n}\rangle=A$.
\end{defn}

\begin{cor}\label{cor:group}
Each group $L$ with $L^{\rm ab}\cong G$ and $[L,L]\cong A$ is isomorphic to $E(\alpha_{\mathbf{a},\mathbf{c}})$ for some admissible pair $(\mathbf{a},\mathbf{c})$.
\end{cor}

\begin{rmk}\label{rmk:alternative}
\rm Let $\mathbf{c}^\vee$ denote the $n\times n$ matrix with entries in $A$ such that the main diagonal entries vanish, the $(i,j)$-entry is $c_{ij}$ if $i<j$, and the $(i,j)$-entry is $-c_{ji}$ if $i>j$.
Then (\ref{eq:eta}) can be rewritten as
\begin{align}
\eta^{\alpha_{\mathbf{a},\mathbf{c}}}(\mathbf{g}^{\mathsf{s}},\mathbf{g}^{\mathsf{t}})=\mathsf{s}^{\rm tr}\mathbf{c}^\vee\mathsf{t}.  \label{eq:eta'}
\end{align}
The multiplication does make sense, despite $A$ being not a ring.

As an immediate consequence of (\ref{eq:commutator}) and (\ref{eq:eta'}),
the center of $E(\alpha_{\mathbf{a},\mathbf{c}})$ is
$$\{(a,\mathbf{g}^{\mathsf{t}})\colon a\in A, \mathsf{t}\in\mathbb{Z}^n, \mathbf{c}^\vee\mathsf{t}=0\}.$$
\end{rmk}

\subsection{Actions of ${\rm End}(G)$ and ${\rm End}(A)$ on $H^2(G,A)$}\label{sec:action}

Let
\begin{align*}
\mathcal{M}_G=\big\{\mathbf{x}\in\mathcal{M}_n(\mathbb{Z}_{p^m})\colon \|\mathbf{x}_{i,j}\|\ge m_i-m_j \ \text{for\ all\ } i>j\big\}.
\end{align*}
By \cite[Theorem 3.3 and Lemma 3.4]{HR07}, there is a ring epimorphism
$$\Theta:\mathcal{M}_G\twoheadrightarrow{\rm End}(G),  \qquad  \mathbf{x}\mapsto\varphi_{\mathbf{x}},$$
where $\varphi_{\mathbf{x}}(\mathbf{g}^{\mathsf{s}})=\mathbf{g}^{\mathbf{x}\mathsf{s}}$ for $\mathsf{s}\in\mathbb{Z}^n$.
The kernel of $\Theta$ consists of $\mathbf{x}$ with $p^{m_i}\mid\mathbf{x}_{i,j}$ for all $i$.

Let $\mathcal{X}_G=\mathcal{M}_G\cap{\rm GL}(n,\mathbb{Z}_{p^m})$.
By \cite[Theorem 3.6]{HR07}, on restriction to $\mathcal{X}_G$, the map $\Theta$ gives rise to a group epimorphism
$\mathcal{X}_G\twoheadrightarrow{\rm Aut}(G)$.

\begin{conv}\label{conv:enumeration}
\rm Fix an enumeration of elements of $\Lambda$, to identify $\Lambda$ with $\{1,\ldots,\check{n}\}$, where $\check{n}={n\choose 2}$.

We shall formally regard elements of $A^n$ as $n$-component row vectors, and regard elements of $A^\Lambda$ as $\check{n}$-component row vectors.
\end{conv}

Given $\mathbf{x}=(x_{ij})\in\mathcal{M}_G$, put
\begin{alignat}{2}
\mathbf{x}^\diamond&=\big(x^\diamond_{(u,v),(i,j)}\big)\in\mathcal{M}_{\check{n}}(\mathbb{Z}_{p^m}), \quad  &&\text{with\ \ \ }
x^\diamond_{(u,v),(i,j)}=x_{ui}x_{vj}-x_{uj}x_{vi};  \label{eq:x-diamond}  \\
\mathbf{x}^\star&=(x^\star_{ij})\in\mathcal{M}_{n}(\mathbb{Z}_{p^m}),  \quad  &&\text{with\ \ \ }
x^\star_{ij}=\begin{cases} [x_{ij}]_{p^{m_i}}/p^{m_i-m_j},&m_i>m_j \\ p^{m_j-m_i}[x_{ij}]_{p^{m_i}},&m_i\le m_j \end{cases}.  \nonumber
\end{alignat}

\begin{rmk}\label{rmk:observe}
\rm 

For any $i,j$, the map ${\rm coker}(\widehat{m_i})\to{\rm coker}(\widehat{m_j})$ induced from $a\mapsto ax^\star_{ij}:=x^\star_{ij}a$ is always well-defined. Hence the map
$${\prod}_{i=1}^n{\rm coker}(\widehat{m_i})\to {\prod}_{i=1}^n{\rm coker}(\widehat{m_i}), \qquad
\overline{\mathbf{a}}\mapsto \overline{\mathbf{a}\mathbf{x}^\star}$$
is a well-defined endomorphism. Here $\overline{\mathbf{a}}$ denotes the image of $\mathbf{a}\in A^n$ under the canonical quotient map $A^n\twoheadrightarrow {\prod}_{i=1}^n{\rm coker}(\widehat{m_i})$.
\end{rmk}

When $p=2$, for $\mathbf{x}\in\mathcal{M}_G$, put $\theta_{\mathbf{x}}=(\theta_{(i,j),h})\in\mathcal{M}_{\check{n}}^n(\mathbb{Z}_{2^m})$, with
$$\theta_{(i,j),h}=\begin{cases}
2^{m_h-1}, &2\nmid x_{ih}x_{jh}\ \text{and}\ m_i=m_j=m_h \\
0,&\text{otherwise} \end{cases}.$$

The following lemma will be useful in proving Theorem \ref{thm:action}:
\begin{lem}\label{lem:assertion}
$\sum\limits_{h=0}^{p^\ell-1}\zeta_{p^\ell}(xh,x)=[x]_{p^\ell}$ for all $x\in\mathbb{Z}_{p^\ell}$.
\end{lem}

\begin{proof}
If $p\nmid x$, then the map of $\{0,1,\ldots,p^\ell-1\}$ to itself defined by $h\mapsto [xh]_{p^\ell}$ is invertible, so
$$\sum\limits_{h=0}^{p^\ell-1}\zeta_{p^\ell}(xh,x)
=\big|\{h\colon 0\le h<p^\ell\ \text{and\ }p^\ell-[x]_{p^\ell}\le[xh]_{p^\ell}<p^\ell\}\big|=[x]_{p^\ell}.$$
In general, write $x=p^ex'$, with $e=\|x\|$ and $p\nmid x'$, and let $\ell'=\ell-e$, then
$$\sum\limits_{h=0}^{p^\ell-1}\zeta_{p^\ell}(xh,x)
=\sum\limits_{h=0}^{p^\ell-1}\zeta_{p^{\ell'}}(x'h,x')=p^e\sum\limits_{h=0}^{p^{\ell'}-1}\zeta_{p^{\ell'}}(x'h,x')
=p^e[x']_{p^{\ell'}}=[x]_{p^\ell}.$$
\end{proof}

\begin{thm}\label{thm:action}
Let $\mathbf{a}\in A^n$ and $\mathbf{c}\in\prod_{(i,j)\in\Lambda}\ker(\widehat{m_i})\subseteq A^\Lambda$.
\begin{enumerate}
  \item[\rm(a)] For each $\psi\in{\rm End}(A)$, denote the endomorphisms of $A^n$ and $A^\Lambda$ induced componentwise from $\psi$ also
        by $\psi$, then $\psi_\ast\alpha_{\mathbf{a},\mathbf{c}}=\alpha_{\psi(\mathbf{a}),\psi(\mathbf{c})}$.
  \item[\rm(b)] For each $\mathbf{x}\in\mathcal{M}_G$,
        $$\varphi_{\mathbf{x}}^\ast[\alpha_{\mathbf{a},\mathbf{c}}]
        =\begin{cases} [\alpha_{\mathbf{a}\mathbf{x}^\star, \mathbf{c}\mathbf{x}^\diamond}], &p>2 \\
        [\alpha_{\mathbf{a}\mathbf{x}^\star+\mathbf{c}\theta_{\mathbf{x}},\mathbf{c}\mathbf{x}^\diamond}],&p=2\end{cases}.$$
\end{enumerate}
\end{thm}

\begin{proof}
The assertion (a) is obvious. We only prove (b).

According to the definition (\ref{eq:alpha}), both sides of the equality are linear in $\mathbf{a},\mathbf{c}$, so it suffices to consider special cases of the following two kinds.

\medskip

(I) Suppose $\mathbf{c}=\mathbf{0}$ and $\mathbf{a}=(a_1,\ldots,a_n)$, with $a_j=\delta_{i,j}a$ for some $1\le i\le n$ and some $a\in A$.

Assume $\varphi_\mathbf{x}^\ast\alpha_{\mathbf{a},\mathbf{0}}=\alpha_{\mathbf{a}',\mathbf{c}'}-\partial\mu$ for some $\mathbf{a}',\mathbf{c}'$ and $\mu:G\to A$.
Since
$$\eta^{\alpha_{\mathbf{a}',\mathbf{c}'}}(\mathbf{g}^{\mathsf{s}},\mathbf{g}^{\mathsf{t}})
=\eta^{\alpha_{\mathsf{a},\mathbf{0}}}(\mathbf{g}^{\mathbf{x}\mathsf{s}},\mathbf{g}^{\mathbf{x}\mathsf{t}})=0$$
for any $\mathsf{s},\mathsf{t}\in\mathbb{Z}^n$, from (\ref{eq:eta}) we see $\mathbf{c}'=\mathbf{0}$.

Write $\mathbf{a}'=(a'_1,\ldots,a'_n)$. Then for any $\mathsf{s},\mathsf{t}$, the following holds in $A$:
\begin{align*}
\mu(\mathbf{g}^{\mathsf{s}+\mathsf{t}})-\mu(\mathbf{g}^{\mathsf{s}})-\mu(\mathbf{g}^{\mathsf{t}})
&=\alpha_{\mathbf{a}',\mathbf{0}}(\mathbf{g}^{\mathsf{s}},\mathbf{g}^{\mathsf{t}})
-(\varphi_\mathbf{x}^\ast\alpha_{\mathbf{a},\mathbf{0}})(\mathbf{g}^{\mathsf{s}},\mathbf{g}^{\mathsf{t}})  \\
&=\sum\limits_{h=1}^n\zeta_{p^{m_h}}(s_h,t_h)a'_h-\zeta_{p^{m_i}}(\mathbf{x}_{i}\mathsf{s},\mathbf{x}_{i}\mathsf{t})a.
\end{align*}
Fix $j$. For $0\le s<p^{m_j}$, taking $\mathsf{s},\mathsf{t}$ with $s_h=s\delta_{h,j}$ and $t_h=\delta_{h,j}$, we obtain
\begin{align*}
\zeta_{p^{m_j}}(s,1)a'_j+\mu(g_j^{s+1})-\mu(g_j^s)-\mu(g_j)=\zeta_{p^{m_i}}(x_{ij}s,x_{ij})a.
\end{align*}
Consequently,
\begin{align}
0&=\sum\limits_{s=0}^{p^{m_j}-1}\big(\mu(g_j^{s+1})-\mu(g_j^s)\big)  \nonumber \\
&=p^{m_j}\mu(g_j)+\sum\limits_{s=0}^{p^{m_j}-1}\zeta_{p^{m_i}}(x_{ij}s,x_{ij})a-\sum\limits_{s=0}^{p^{m_j}-1}\zeta_{p^{m_j}}(s,1)a'_j.
\label{eq:important}
\end{align}

If $i>j$, then $\zeta_{p^{m_i}}(x_{ij}s,x_{ij})=\zeta_{p^{m_j}}(x^\star_{ij}s,x^\star_{ij})$,
so by Lemma \ref{lem:assertion},
$$\sum_{s=0}^{p^{m_j}-1}\zeta_{p^{m_i}}(x_{ij}s,x_{ij})=\sum_{s=0}^{p^{m_j}-1}\zeta_{p^{m_j}}(x^\star_{ij}s,x^\star_{ij})
=[x^\star_{ij}]_{p^{m_j}}=x^\star_{ij};$$
if $i\le j$, then since $\zeta_{p^{m_i}}(x_{ij}s,x_{ij})$ is periodic in $s$ with period $p^{m_i}$, we have
$$\sum_{s=0}^{p^{m_j}-1}\zeta_{p^{m_i}}(x_{ij}s,x_{ij})=p^{m_j-m_i}\sum_{s=0}^{p^{m_i}-1}\zeta_{p^{m_i}}(x_{ij}s,x_{ij})
=p^{m_j-m_i}[x_{ij}]_{p^{m_i}}=x^\star_{ij}.$$

Hence (\ref{eq:important}) implies $a'_j=p^{m_j}\mu(g_j)+x^\star_{ij}a$.

We have shown $a'_j-x^\star_{ij}a\in\im(\widehat{m_j})$ for each $j$,
so that $\mathbf{a}'-\mathbf{a}\mathbf{x}^\star\in\prod_{j=1}^n\im(\widehat{m_j})$.
Thus by Theorem \ref{thm:H2},
$$\varphi_\mathbf{x}^\ast[\alpha_{\mathbf{a},\mathbf{0}}]=[\alpha_{\mathbf{a}',\mathbf{0}}]
=[\alpha_{\mathbf{a}\mathbf{x}^\star,\mathbf{0}}].$$

\medskip

(II) Suppose $\mathbf{a}=\mathbf{0}$ and $\mathbf{c}=\{c_{ij}\}_{(i,j)\in\Lambda}$, with $c_{ij}=\delta_{(i,j),(u,v)}c$ for some $(u,v)\in\Lambda$ and some $c\in\ker(\hat{u})$.

Assume $\varphi_\mathbf{x}^\ast\alpha_{\mathbf{0},\mathbf{c}}=\alpha_{\mathbf{a}',\mathbf{c}'}-\partial\mu$
for some map $\mu:G\to A$ and some $\mathbf{a}'=(a'_1,\ldots,a'_n)$, $\mathbf{c}'=\{c'_{ij}\}_{(i,j)\in\Lambda}$.
Then for any $\mathsf{s},\mathsf{t}\in\mathbb{Z}^n$,
\begin{align}
\mu(\mathbf{g}^{\mathsf{s}+\mathsf{t}})-\mu(\mathbf{g}^{\mathsf{s}})-\mu(\mathbf{g}^{\mathsf{t}})
=\sum\limits_{i=1}^n\zeta_{p^{m_i}}(s_i,t_i)a'_i-\sum\limits_{1\le i<j\le n}t_is_jc'_{ij}
-(\mathbf{x}_{u}\mathsf{t})(\mathbf{x}_{v}\mathsf{s})c.  \label{eq:difference-II}
\end{align}
Interchanging $\mathsf{s}$ with $\mathsf{t}$ and taking the difference between the resulting equation and (\ref{eq:difference-II}), we arrive at
\begin{align*}
\sum\limits_{1\le i<j\le n}(s_it_j-t_is_j)c'_{ij}
=\sum\limits_{1\le i<j\le n}(x_{ui}x_{vj}-x_{uj}x_{vi})(s_it_j-t_is_j)c.
\end{align*}
Hence for each $(i,j)\in\Lambda$, we have
$$c'_{ij}=(x_{ui}x_{vj}-x_{uj}x_{vi})c=x^\diamond_{(u,v),(i,j)}\cdot c,$$
which is exactly the entry of $\mathbf{c}\mathbf{x}^\diamond$ indexed by $(i,j)$. This shows $\mathbf{c}'=\mathbf{c}\mathbf{x}^\diamond$.

Fix $h$. For $0\le s<p^{m_h}$, setting $s_j=s\delta_{j,h}$ and $t_j=\delta_{j,h}$ in (\ref{eq:difference-II}), we obtain
$$\mu(g_h^{s+1})-\mu(g_h^s)=\mu(g_h)+\zeta_{p^{m_h}}(s,1)a'_h+sx_{uh}x_{vh}c.$$
Consequently,
\begin{align*}
0&=\sum_{s=0}^{p^{m_h}-1}\big(\mu(g_h^{s+1})-\mu(g_h^s)\big)  \\
&=p^{m_h}\mu(g_h)+\sum_{s=0}^{p^{m_h}-1}\zeta_{p^{m_h}}(s,1)a'_h+{p^{m_h}\choose 2}x_{uh}x_{vh}c \\
&=a'_h+p^{m_h}\mu(g_h)+\begin{cases} p^{m_h}\cdot\frac{1}{2}(p^{m_h}-1)x_{uh}x_{vh}c, &p>2 \\
2^{m_h}\cdot 2^{m_h-1}x_{uh}x_{vh}c-2^{m_h-1}x_{uh}x_{vh}c,&p=2 \end{cases}.
\end{align*}
It follows that $a'_h\in\im(\widehat{m_h})$ when $p>2$.

When $p=2$, we have
$$2^{m_h-1}x_{uh}x_{vh}=\theta_{(u,v),h}.$$
Indeed, if $2^{m_h-1}x_{uh}x_{vh}c\ne 0$, then $|c|\ge 2^{m_h+\|x_{uh}\|+\|x_{vh}\|}$; now that $c\in\ker(\hat{u})$, i.e. $2^{m_u}c=0$, one has
$$m_u\ge m_h+\|x_{uh}\|+\|x_{vh}\|\ge m_h+\max\{0,m_u-m_h\}+\max\{0,m_v-m_h\},$$
which together with $m_v\ge m_u$ force $m_u=m_v=m_h$ and $2\nmid x_{uh}x_{vh}$.

Thus $\varphi_\mathbf{x}^\ast[\alpha_{\mathbf{0},\mathbf{c}}]
=[\alpha_{\mathbf{a}',\mathbf{c}'}]=[\alpha_{\mathbf{c}\theta_{\mathbf{x}},\mathbf{c}\mathbf{x}^\diamond}].$
\end{proof}

\section{Homomorphisms between central extensions}

In this section, assume $p>2$. Write elements of $A=\mathbb{Z}_{p^{d_1}}\times\cdots\times \mathbb{Z}_{p^{d_r}}$ in the forms of $r$-component column vectors. Continue to adopt Convention \ref{conv:enumeration}.
For $(\mathbf{a},\mathbf{c})\in A^n\times A^\Lambda$, write $\mathbf{a},\mathbf{c}$ respectively as a $r\times n$, $r\times {n\choose 2}$ matrix.

Recall $d=d_r$. Let
$$\mathcal{M}^A=\big\{\mathbf{y}\in\mathcal{M}_r(\mathbb{Z}_{p^d})\colon \|\mathbf{y}_{i,j}\|\ge d_j-d_i\ \text{for\ all\ }j>i\big\}.$$
Similarly as \cite[Theorem 3.3]{HR07}, each $\mathbf{y}\in\mathcal{M}^A$ defines $\psi^{\mathbf{y}}\in{\rm End}(A)$ by $a\mapsto\mathbf{y}a$, and conversely, each endomorphism of $A$ has this form.

\subsection{Lifting endomorphisms}

Suppose $(\mathbf{a},\mathbf{c})$, $(\mathbf{a}',\mathbf{c}')$ are admissible pairs.

Given $\mathbf{x}\in\mathcal{M}_G$, by Proposition \ref{prop:basic} (c), $\varphi_{\mathbf{x}}\in{\rm End}(G)$ can be lifted to
$E(\alpha_{\mathbf{a},\mathbf{c}})\to E(\alpha_{\mathbf{a}',\mathbf{c}'})$
if and only if there exists $\mathbf{y}\in\mathcal{M}^A$ with $$\varphi_{\mathbf{x}}^\ast[\alpha_{\mathbf{a}',\mathbf{c}'}]=\psi^{\mathbf{y}}_\ast[\alpha_{\mathbf{a},\mathbf{c}}],$$
which, by Theorem \ref{thm:action}, is equivalent to
\begin{align}
\overline{\mathbf{a}'\mathbf{x}^\star}=\overline{\mathbf{y}\mathbf{a}},  \qquad
\mathbf{c}'\mathbf{x}^\diamond=\mathbf{y}\mathbf{c}.   \label{eq:compatible}
\end{align}
Recall that by Remark \ref{rmk:observe}, $\overline{\mathbf{a}}$ denotes the image of $\mathbf{a}$ under the quotient map
$A^n\twoheadrightarrow\prod_{i=1}^n{\rm coker}(\widehat{m_i})$.
By Remark \ref{rmk:determine}, $\mathbf{y}$ (if exsits) is determined by $\mathbf{x}$.

Focusing on automorphisms, let
\begin{align}
\mathcal{A}_{\mathbf{a},\mathbf{c}}&=\{\mathbf{x}\in\mathcal{X}_G\colon \overline{\mathbf{a}\mathbf{x}^\star}=\overline{\mathbf{y}\mathbf{a}}, \mathbf{c}\mathbf{x}^\diamond=\mathbf{y}\mathbf{c}\ \text{for\ some\ }\mathbf{y}\}, \label{eq:def-A}  \\ \mathcal{L}_{\mathbf{a},\mathbf{c}}&=\{\varphi_{\mathbf{x}}\colon \mathbf{x}\in\mathcal{A}_{\mathbf{a},\mathbf{c}}\}.  \nonumber
\end{align}

\begin{thm}\label{thm:auto-group}
There is a short exact sequence
$$1\to\hom(G,A)\stackrel{c}\longrightarrow {\rm Aut}(E(\alpha_{\mathbf{a},\mathbf{c}}))\stackrel{q}\longrightarrow
\mathcal{L}_{\mathbf{a},\mathbf{c}}\to 1,$$
where $c$ sends $\sigma:G\to A$ to the automorphism $(a,z)\mapsto(a+\sigma(z),z)$, and
$q$ sends an automorphism of $E(\alpha_{\mathbf{a},\mathbf{c}})$ to its induced automorphism on $G$.
\end{thm}

This was already known (see \cite[Page 1]{Ro77}). Here we give a proof, for the reader's convenience.

\begin{proof}
Clearly $c$ is injective, $q$ is surjective, and $q\circ c$ is trivial. It suffices to show $\ker(q)\subseteq\im(c)$.

Suppose $\phi\in \ker(q)$, i.e. $\phi$ is an automorphism inducing the identity on $G$. By Remark \ref{rmk:determine}, $\phi|_A={\rm id}$.
By Proposition \ref{prop:basic} (c), there exists a map $\mu:G\to A$ with $\partial\mu=0$ such that $\phi(a,z)=(a+\mu(z),z)$. The condition $\partial\mu=0$ is equivalent to $\mu\in\hom(G,A)$. Thus, $\phi=c(\mu)$.
\end{proof}

\subsection{An efficient method for computation}

Define an embedding $\iota:A\hookrightarrow \mathbb{Z}_{p^{d}}^r$ by
\begin{align}
\big(w_1,\ldots,w_r)^{\rm tr}\mapsto(p^{d-d_1}w_1,\ldots,p^{d-d_r}w_r\big)^{\rm tr},  \label{eq:embedding-1}
\end{align}
under which $A$ can be identified with the subgroup $\widetilde{A}$ of $\mathbb{Z}_{p^d}^r$, with
\begin{align*}
\widetilde{A}=\{(u_1,\ldots,u_r)^{\rm tr}\colon u_j\in\mathbb{Z}_{p^d}, \|u_j\|\ge d-d_j\}.
\end{align*}

For each $1\le i\le n$, the image of $\ker(\widehat{m_i})$ under (\ref{eq:embedding-1}) is
\begin{align*}
Q_i:=\{(v_1,\ldots,v_r)^{\rm tr}\in \mathbb{Z}_{p^d}^r\colon\|v_j\|\ge d-\min\{m_i,d_j\}, \ 1\le j\le r \}.
\end{align*}
Define $\pi_i:A\to\mathbb{Z}_{p^d}^r$ by
\begin{align*}
(w_1,\ldots,w_r)^{\rm tr}\mapsto\big(p^{d-\min\{d_1,m_i\}}w_1,\ldots,p^{d-\min\{d_r,m_i\}}w_r\big)^{\rm tr}.
\end{align*}
Clearly, $\im(\pi_i)=Q_i$, and $\ker(\pi_i)=p^{m_i}A$, so $\pi_i$ induces an embedding
${\rm coker}(\widehat{m_i})\hookrightarrow\mathbb{Z}_{p^d}^r$, whose image is $Q_i$.

\begin{rmk}
\rm Both $\ker(\widehat{m_i})$ and ${\rm coker}(\widehat{m_i})$ can be identified with $Q_i\leqslant\widetilde{A}$.

Working with $Q_i$ in place of ${\rm coker}(\widehat{m_i})$, we get rid of the indeterminacy caused by the representative of an element of
${\rm coker}(\widehat{m_i})$. 
\end{rmk}

Given an admissible pair $(\mathbf{a},\mathbf{c})$, say $\mathbf{a}=(a_1,\ldots,a_n)$, $\mathbf{c}=\{c_{ij}\}_{(i,j)\in\Lambda}$, let $$\pi(\mathbf{a})=(\pi_1(a_1),\ldots,\pi_n(a_n)), \qquad \iota(\mathbf{c})=\{\iota(c_{ij})\}_{(i,j)\in\Lambda}.$$
Due to (\ref{eq:generate}), we have $\langle\iota(c_{12}),\ldots,\iota(c_{n-1,n})\rangle=\widetilde{A}$, so $\pi(\mathbf{a})=\iota(\mathbf{c})\mathbf{w}$ for some $\mathbf{w}\in\mathcal{M}_{\check{n}}^n(\mathbb{Z}_{p^d})$, with $\check{n}={n\choose 2}$.
Let
$\mathsf{N}=\{\mathsf{u}\in \mathbb{Z}_{p^d}^\Lambda\colon\iota(\mathbf{c})\mathsf{u}=0\}.$

Given another admissible pair $(\mathbf{a}',\mathbf{c}')$, write $\pi(\mathbf{a}')=\iota(\mathbf{c}')\mathbf{w}'$ for some $\mathbf{w}'\in\mathcal{M}_{\check{n}}^n(\mathbb{Z}_{p^d})$, and let
$\mathsf{N}'=\{\mathsf{u}\in \mathbb{Z}_{p^d}^\Lambda\colon\iota(\mathbf{c}')\mathsf{u}=0\}.$

The condition (\ref{eq:compatible}), which is equivalent to that $\varphi_{\mathbf{x}}$ can be lifted, becomes
\begin{align}
\pi(\mathbf{a}')\mathbf{x}^\star&=\mathbf{y}\pi(\mathbf{a}),   \label{eq:compatible1} \\
\iota(\mathbf{c}')\mathbf{x}^\diamond&=\mathbf{y}\iota(\mathbf{c}). \label{eq:compatible2}
\end{align}

\begin{lem}\label{lem:criterion}
Given $\mathbf{x}\in\mathcal{M}_G$, the endomorphism $\varphi_{\mathbf{x}}\in{\rm End}(G)$ can be lifted to a homomorphism $E(\alpha_{\mathbf{a},\mathbf{c}})\to E(\alpha_{\mathbf{a}',\mathbf{c}'})$ if and only if $\mathbf{x}^\diamond\mathsf{N}\leqslant\mathsf{N}'$ and $\iota(\mathbf{c}')(\mathbf{w}'\mathbf{x}^\star-\mathbf{x}^\diamond\mathbf{w})=0$.
If $\varphi_{\mathbf{x}}\in{\rm Aut}(G)$, then each lifting is an isomorphism.
\end{lem}

\begin{proof}
Suppose (\ref{eq:compatible1}) and (\ref{eq:compatible2}) hold.
Then $\iota(\mathbf{c}')\mathbf{x}^\diamond\mathsf{N}=\mathbf{y}\iota(\mathbf{c})\mathsf{N}=0$, implying $\mathbf{x}^\diamond\mathsf{N}\leqslant\mathsf{N}'$,
and
$$\iota(\mathbf{c}')(\mathbf{w}'\mathbf{x}^\star-\mathbf{x}^\diamond\mathbf{w})
=\iota(\mathbf{c}')\mathbf{w}'\mathbf{x}^\star-\mathbf{y}\iota(\mathbf{c})\mathbf{w}=\pi(\mathbf{a}')\mathbf{x}^\star-\mathbf{y}\pi(\mathbf{a})=0.$$

Conversely, suppose $\mathbf{x}^\diamond\mathsf{N}\leqslant\mathsf{N}'$ and $\iota(\mathbf{c}')(\mathbf{w}'\mathbf{x}^\star-\mathbf{x}^\diamond\mathbf{w})=0$.
The first condition ensures that the map
$\iota(\mathbf{c})\mathsf{u}\mapsto\iota(\mathbf{c}')\mathbf{x}^\diamond\mathsf{u}$
is well-defined, so there exists $\mathbf{y}\in\mathcal{M}^A$ such that $\mathbf{y}\iota(\mathbf{c})\mathsf{u}=\iota(\mathbf{c}')\mathbf{x}^\diamond\mathsf{u}$ for all $\mathsf{u}\in \mathbb{Z}_{p^d}^\Lambda$.
Hence (\ref{eq:compatible2}) holds, and then (\ref{eq:compatible1}) follows as
$$\mathbf{y}\pi(\mathbf{a})=\mathbf{y}\iota(\mathbf{c})\mathbf{w}=\iota(\mathbf{c}')\mathbf{x}^\diamond\mathbf{w}
=\iota(\mathbf{c}')\mathbf{w}'\mathbf{x}^\star=\pi(\mathbf{a}')\mathbf{x}^\star.$$

Now suppose $\varphi_\mathbf{x}$ is invertible,
and $\phi:E(\alpha_{\mathbf{a},\mathbf{c}})\to E(\alpha_{\mathbf{a}',\mathbf{c}'})$ is a lifting of $\varphi_\mathbf{x}$.
Observe that the abelianization $E(\alpha_{\mathbf{a}',\mathbf{c}'})\twoheadrightarrow G$ takes $\im(\phi)$ to
$\im(\varphi_\mathbf{x})=G$. By Proposition \ref{prop:basic} (a), $\im(\phi)=E(\alpha_{\mathbf{a}',\mathbf{c}'})$.
So actually $\phi$ is invertible.
\end{proof}

In practical computation, we may choose a generating set $\{\mathsf{u}_1,\ldots,\mathsf{u}_\ell\}$ of $\mathsf{N}$, and express the condition $\mathbf{x}^\diamond\mathsf{N}\leqslant\mathsf{N}'$ as
$\iota(\mathbf{c}')\mathbf{x}^\diamond\mathsf{u}_i=0$, $1\le i\le\ell.$

In particular,
\begin{align*}
\mathcal{A}_{\mathbf{a},\mathbf{c}}=\big\{\mathbf{x}\in\mathcal{X}_G\colon \iota(\mathbf{c})\mathbf{x}^\diamond\mathsf{u}_1=\cdots=\iota(\mathbf{c})\mathbf{x}^\diamond\mathsf{u}_\ell=0, \
\iota(\mathbf{c})(\mathbf{w}\mathbf{x}^\star-\mathbf{x}^\diamond \mathbf{w})=0\big\}.
\end{align*}
This characterizes $\mathbf{x}\in\mathcal{A}_{\mathbf{a},\mathbf{c}}$ without mentioning $\mathbf{y}$.

\begin{exmp}
\rm Let $A=\mathbb{Z}_{p}\times \mathbb{Z}_{p^3}$, $G=\mathbb{Z}_{p^2}\times \mathbb{Z}_{p^3}\times \mathbb{Z}_{p^4}$.
As is easy to see,
\begin{align*}
Q_1&=\{(v_1,v_2)^{\rm tr}\in \mathbb{Z}_{p^3}^2\colon \|v_1\|\ge 2, \|v_2\|\ge 1\},   \\
Q_2&=Q_3=\{(v_1,v_2)^{\rm tr}\in \mathbb{Z}_{p^3}^2\colon \|v_1\|\ge 2\}.
\end{align*}

Enumerate the elements of $\Lambda$ as $(2,3)$, $(1,3)$, $(1,2)$.
Let
\begin{alignat*}{2}
\mathbf{a}&=\left(\begin{array}{ccc} 1 & 0 & 3 \\ 0 & 1 & p^2 \end{array}\right), \qquad
&\mathbf{c}&=\left(\begin{array}{ccc} 0 & 1 & 0 \\ 1 & p & -p \end{array}\right);  \\
\mathbf{a}'&=\left(\begin{array}{ccc} 0 & 5 & 0 \\ 1 & p & 2 \end{array}\right), \qquad
&\mathbf{c}'&=\left(\begin{array}{ccc} 4 & 0 & 1 \\ 1 & 0 & 0 \end{array}\right).
\end{alignat*}
We emphasize that the first, second, third columns of $\mathbf{c}$ are respectively $c_{23}$, $c_{13}$, $c_{12}$, if $\mathbf{c}=\{c_{ij}\}_{(i,j)\in\Lambda}$.
Note that $\pi_1$ sends $(w_1,w_2)^{\rm tr}$ to $(p^2w_1,pw_2)^{\rm tr}$, and $\iota,\pi_2,\pi_3$ send
$(w_1,w_2)^{\rm tr}$ to $(p^2w_1,w_2)^{\rm tr}$. Hence
\begin{alignat*}{2}
\pi(\mathbf{a})&=\left(\begin{array}{ccc} p^2 & 0 & 3p^2 \\ 0 & 1 & p^2 \end{array}\right), \qquad
&&\iota(\mathbf{c})=\left(\begin{array}{ccc} 0 & p^2 & 0 \\ 1 & p & -p \end{array}\right);  \\
\pi(\mathbf{a}')&=\left(\begin{array}{ccc} 0 & 5p^2 & 0 \\ p & p & 2 \end{array}\right), \qquad
&&\iota(\mathbf{c}')=\left(\begin{array}{ccc} 4p^2 & 0 & p^2 \\ 1 & 0 & 0 \end{array}\right).
\end{alignat*}
We may take
\begin{alignat*}{2}
\mathbf{w}=\left(\begin{array}{ccc} -p & 1 & p^2-3p \\ 1 & 0 & 3  \\  0 & 0 & 0 \end{array}\right), \qquad
\mathbf{w}'=\left(\begin{array}{ccc} p & p & 2 \\  0 & 0 &  0  \\ -4p & 5-4p & -8 \end{array}\right),
\end{alignat*}
and take $\mathsf{u}_1=(p,0,1)^{\rm tr}$, $\mathsf{u}_2=(p^2,-p,0)^{\rm tr}$ as generators for $\mathsf{N}$.

Each $\mathbf{x}=(x_{ij})\in\mathcal{M}_G$ satisfies $p\mid x_{21},x_{32}$ and $p^2\mid x_{31}$.
We have
\begin{align*}
\mathbf{x}^\star&=\left(\begin{array}{ccc} x_{11} & px_{12} & p^2x_{13} \\ {[x_{21}]_{p^3}}/p & x_{22} & px_{23}  \\
{[x_{31}]_{p^4}}/p^2& {[x_{32}]_{p^4}}/p & x_{33} \end{array}\right), \\
\mathbf{x}^\diamond&=\left(\begin{array}{ccc} x_{22}x_{33}-x_{23}x_{32} & x_{21}x_{33}-x_{23}x_{31} & x_{21}x_{32}-x_{22}x_{31} \\
x_{12}x_{33}-x_{13}x_{32} & x_{11}x_{33}-x_{13}x_{31} & x_{11}x_{32}-x_{12}x_{31}  \\
x_{12}x_{23}-x_{13}x_{22} & x_{11}x_{23}-x_{13}x_{21} & x_{11}x_{22}-x_{12}x_{21} \end{array}\right).
\end{align*}
Then the conditions
$\iota(\mathbf{c}')\mathbf{x}^\diamond\mathsf{u}_1=\iota(\mathbf{c}')\mathbf{x}^\diamond\mathsf{u}_2=0$ and
$\iota(\mathbf{c}')(\mathbf{w}'\mathbf{x}^\star-\mathbf{x}^\diamond\mathbf{w})=0$ can be explicitly written down as equations in $\mathbb{Z}_{p^3}$, which can be used to decide whether $\varphi_{\mathbf{x}}$ is liftable.
\end{exmp}

\subsection{Towards classifying $p$-groups of class $2$}

Let $\mathcal{E}(G,A)$ denote the set of $p$-groups $L$ such that $[L,L]\cong A$ and $L^{\rm ab}\cong G$.

Define an equivalence relation $\asymp$ among admissible pairs by declaring $(\mathbf{a},\mathbf{c})\asymp(\mathbf{a}',\mathbf{c}')$ if $\overline{\mathbf{y}\mathbf{a}}=\overline{\mathbf{a}'\mathbf{x}^\star}$, $\mathbf{y}\mathbf{c}=\mathbf{c}'\mathbf{x}^\diamond$ for some
$(\mathbf{x},\mathbf{y})\in\mathcal{X}_G\times\mathcal{X}^A$, where
$$\mathcal{X}^A=\big\{\mathbf{y}\in{\rm GL}(r,\mathbb{Z}_{p^d})\colon \|\mathbf{y}_{i,j}\|\ge d_j-d_i\ \text{for\ all\ }j>i\big\}.$$
As we have shown, to classify the groups in $\mathcal{E}(G,A)$ up to isomorphism, it suffices to classify admissible pairs up to equivalence.

An alternative approach will be more convenient. For
$$\tilde{\mathbf{a}}=(\tilde{a}_1,\ldots,\tilde{a}_n)\in\widetilde{A}^n, \qquad \tilde{\mathbf{c}}=\{\tilde{c}_{ij}\}_{(i,j)\in\Lambda}\in\widetilde{A}^\Lambda,$$
we call $(\tilde{\mathbf{a}},\tilde{\mathbf{c}})$ an {\it admissible pair} if $\tilde{a}_i,\tilde{c}_{ij}\in Q_i$ and
$\langle \tilde{c}_{12},\ldots,\tilde{c}_{n-1,n}\rangle=\widetilde{A}$.
Declare $(\tilde{\mathbf{a}},\tilde{\mathbf{c}})\asymp(\tilde{\mathbf{a}}',\tilde{\mathbf{c}}')$ if $\mathbf{y}\tilde{\mathbf{a}}=\tilde{\mathbf{a}}'\mathbf{x}^\star$, $\mathbf{y}\tilde{\mathbf{c}}=\tilde{\mathbf{c}}'\mathbf{x}^\diamond$ for some
$(\mathbf{x},\mathbf{y})\in\mathcal{X}_G\times\mathcal{X}^A$.
In view of (\ref{eq:compatible1}), (\ref{eq:compatible2}), the classification problem can also be settled by finding a complete list of representatives of equivalence classes of admissible pairs in $\widetilde{A}^n\times\widetilde{A}^\Lambda$.


The classification procedure can be implemented in two steps.
\begin{enumerate}
  \item Define an equivalence relation $\sim$ among all $\tilde{\mathbf{c}}=\{\tilde{c}_{ij}\}_{(i,j)\in\Lambda}\in\widetilde{A}^\Lambda$ with $\tilde{c}_{ij}\in Q_i$ and $\langle \tilde{c}_{12},\ldots,\tilde{c}_{n-1,n}\rangle=\widetilde{A}$ by declaring $\tilde{\mathbf{c}}\sim\tilde{\mathbf{c}}'$ if $\tilde{\mathbf{c}}'=\mathbf{y}^{-1}\tilde{\mathbf{c}}\mathbf{x}^\diamond$ for some
      $(\mathbf{x},\mathbf{y})\in\mathcal{X}_G\times\mathcal{X}^A$. Find a ``normal" representative (whose sense depends on the situation) for each equivalence class.
  \item Fix a normal $\tilde{\mathbf{c}}$. Define an equivalence relation $\approx$ among all $\tilde{\mathbf{a}}=(\tilde{a}_1,\ldots,\tilde{a}_n)\in\widetilde{A}^n$ with $\tilde{a}_i\in Q_i$ by declaring $\tilde{\mathbf{a}}\approx\tilde{\mathbf{a}}'$ if $\tilde{\mathbf{a}}'=\mathbf{y}\tilde{\mathbf{a}}(\mathbf{x}^\star)^{-1}$, for some $(\mathbf{x},\mathbf{y})\in\mathcal{X}_G\times\mathcal{X}^A$ with $\mathbf{y}^{-1}\tilde{\mathbf{c}}\mathbf{x}^\diamond=\tilde{\mathbf{c}}$.
      Find a ``normal" representative for each equivalence class.
\end{enumerate}
When these finish, for each admissible pair $(\tilde{\mathbf{a}},\tilde{\mathbf{c}})$ such that $\tilde{\mathbf{a}}$, $\tilde{\mathbf{c}}$ are normal, take $(\mathbf{a},\mathbf{c})\in A^n\times A^\Lambda$ with $\pi(\mathbf{a})=\tilde{\mathbf{a}}$ and $\iota(\mathbf{c})=\tilde{\mathbf{c}}$.
Then each group in $\mathcal{E}(G,A)$ is isomorphic to exactly one $E(\alpha_{\mathbf{a},\mathbf{c}})$ for such $(\mathbf{a},\mathbf{c})$.

Section 4.1 is an illustration of this procedure.

\section{Two-generator $p$-groups of class $2$ and the orders of their automorphism groups, for $p>2$}

\subsection{Recovering the classification}

By Remark \ref{rmk:invariants}, we may suppose $A=\mathbb{Z}_{p^d}$ and $G=\mathbb{Z}_{p^{m_1}}\times \mathbb{Z}_{p^{m_2}}$, with $d\le m_1\le m_2$.
Then $|\Lambda|=1$, $Q_1=Q_2=\widetilde{A}={\rm coker}(\widehat{m_1})={\rm coker}(\widehat{m_2})=A$.

For an admissible pair $(\mathbf{a},\mathbf{c})\in A^2\times A$, we can write $\mathbf{a}=(a_1,a_2)$ and $\mathbf{c}=c$, with $a_1,a_2,c\in \mathbb{Z}_{p^d}$. The condition (\ref{eq:generate}) is simply equivalent to $p\nmid c$.

For $(\mathbf{x},\mathbf{y})\in\mathcal{X}_G\times\mathcal{X}^A$, we can write $\mathbf{x}=(x_{ij})$, with $p\nmid \det(\mathbf{x})$ and $\|x_{21}\|\ge m_2-m_1$, and write $\mathbf{y}=y\in \mathbb{Z}_{p^d}$ with $p\nmid y$.
Note that $\mathbf{x}^\diamond=\det(\mathbf{x})$.

Up to the equivalence $\sim$, we may assume $c=1$.

Fix $\mathbf{c}=c=1$. 
The condition $\mathbf{y}^{-1}\mathbf{c}\mathbf{x}^\diamond=\mathbf{c}$ is equivalent to $y=\det(\mathbf{x})$. Hence, for any $\mathbf{a},\mathbf{a}'\in A^2$, $\mathbf{a}\approx\mathbf{a}'$ if and only if there exists $\mathbf{x}\in\mathcal{X}_G$ such that
\begin{align*}
\mathbf{a}'=y\mathbf{a}(\mathbf{x}^\star)^{-1}
=\mathbf{a}\left(\begin{array}{cc} x_{22} & -p^{m_2-m_1}x_{12} \\ -[x_{21}]_{p^{m_2}}/p^{m_2-m_1} & x_{11}\end{array}\right).  
\end{align*}

If $m_1=m_2$, then $(a_1,a_2)\approx(0,p^{\min\{\|a_1\|,\|a_2\|\}})$, which can be written as $(p^d,p^k)$, with $0\le k\le d$.
When $k\ne k'$, it is easy to see $(0,p^k)\not\approx(0,p^{k'})$. 

Now suppose $m_2>m_1$. Then $(a_1,a_2)\approx(p^{\|a_1\|},p^{\|a_2\|})$, which can be written as $(p^{k_1},p^{k_2})$, with $0\le k_1,k_2\le d$.
Furthermore,
\begin{itemize}
  \item if $k_2\le k_1$, then $(p^{k_1},p^{k_2})\approx(0,p^{k_2})$;
  \item if $d>k_2\ge k_1+m_2-m_1$, then $(p^{k_1},p^{k_2})\approx(p^{k_1},0)$.
\end{itemize}
Thus, each $\mathbf{a}\in A^2$ is equivalent to $(p^{k_1},p^{k_2})$ for $(k_1,k_2)$ satisfying one of the following {\it reduced} conditions:
\begin{enumerate}
  \item[\rm(a)] $0\le k_2\le k_1=d$;
  \item[\rm(b)] $0\le k_1<k_2=d$;
  \item[\rm(c)] $0\le k_1<k_2<\min\{d,k_1+m_2-m_1\}$.
\end{enumerate}
We show that the $(p^{k_1},p^{k_2})$'s for distinct pairs $(k_1,k_2)$ are non-equivalent to each another.
Assume $(k_1,k_2)$, $(k'_1,k'_2)$ are reduced and
\begin{align*}
(p^{k'_1},p^{k'_2})&=(p^{k_1},p^{k_2})\left(\begin{array}{cc} x_{22} & -p^{m_2-m_1}x_{12} \\ -[x_{21}]_{p^{m_2}}/p^{m_2-m_1} & x_{11}\end{array}\right)    \\
&=\big(p^{k_1}x_{22}+p^{k_2}z,p^{k_2}x_{11}-p^{k_1+m_2-m_1}x_{12}\big),   \qquad z=-\frac{[x_{21}]_{p^{m_2}}}{p^{m_2-m_1}}.
\end{align*}
\begin{enumerate}
  \item If $0\le k_2\le k_1=d$, then $k'_2=\|p^{k_2}x_{11}\|=k_2$, and $k'_1=\|p^{k_2}z\|\ge k_2$; since $(k'_1,k'_2)$ is reduced,
        we have $k'_1=d$.
  \item If $0\le k_1<k_2=d$, then $k'_1=\|p^{k_1}x_{22}\|=k_1$, $k'_2=\|p^{k_1+m_2-m_1}x_{12}\|\ge k_1+m_2-m_1$;
        since $(k'_1,k'_2)$ is reduced, we have $k'_2=d$.
  \item If $0\le k_1<k_2<\min\{d,k_1+m_2-m_1\}$, then $k'_1=\|p^{k_1}x_{22}+p^{k_2}z\|=k_1$,
        and $k'_2=\|p^{k_2}x_{11}-p^{k_1+m_2-m_1}x_{12}\|=k_2$.
\end{enumerate}

Let $H_{k_1,k_2}$ denote $E(\alpha_{(p^{k_1},p^{k_2}),1})$. We have established
\begin{thm}\label{thm:classification}
Suppose $L$ is a $p$-group of class 2 such that $[L,L]\cong \mathbb{Z}_{p^d}$ and $L^{\rm ab}\cong \mathbb{Z}_{p^{m_1}}\times\mathbb{Z}_{p^{m_2}}$, with $d\le m_1\le m_2$.
\begin{enumerate}
  \item[\rm 1.] If $m_1=m_2=m$, then $L\cong H_{d,k}$ for a unique $k$ with $0\le k\le d$.
  \item[\rm 2.] If $m_1<m_2$, then $L\cong H_{k_1,k_2}$ for a unique pair $(k_1,k_2)$ satisfying one of the following conditions:
  \begin{enumerate}
    \item[\rm(a)] $0\le k_2\le k_1=d$;
    \item[\rm(b)] $0\le k_1<k_2=d$;
    \item[\rm(c)] $0\le k_1<k_2<\min\{d,k_1+m_2-m_1\}$.
  \end{enumerate}
\end{enumerate}
\end{thm}

\bigskip

The commutator used in \cite{AMM12} was defined by $[a,b]=a^{-1}b^{-1}ab$.
In our notation of commutator, the group denoted by $(m_2,m_1,d;k_2,k_1)$ in \cite[Theorem 1.1]{AMM12} can be presented by (with $u=a^{-1},v=b^{-1}$)
$$\big\langle u,v\mid [u,v]^{p^d}=[u,v,u]=[u,v,u]=1,\ u^{-p^{m_2}}=[u,v]^{p^{k_2}},\ v^{-p^{m_1}}=[u,v]^{p^{k_1}}\big\rangle.$$

\begin{lem}
The homomorphism $\vartheta:(m_2,m_1,d;k_2,k_1)\to H_{k_1,k_2}$ determined by $u\mapsto(0,g_2)$, $v\mapsto(0,g_1)$ is an isomorphism.
\end{lem}

\begin{proof}
Let $\alpha=\alpha_{(p^{k_1},p^{k_2}),1}$. By (\ref{eq:eta}), $\eta^{\alpha}(g_2,g_1)=-1$, so
$$\vartheta([u,v])=[(0,g_2),(0,g_1)]=(\eta^{\alpha}(g_2,g_1),1)=(-1,1).$$
Using that $\alpha$ is normalized, it is easy to compute $(-1,1)^{\ell}=(-\ell,1)$, for all $\ell\in\mathbb{Z}_{>0}$.
In particular, $\vartheta([u,v])^{p^d}=(-p^d,1)=(0,1)$.

For $i=1,2$ and each $0\le h<p^{m_i}$,
$$\alpha(g_i,g_i^h)=\zeta_{p^{m_i}}(1,h)\cdot p^{k_i}=\delta_{h,p^{m_i}-1}\cdot p^{k_i}.$$
Hence
$$(0,g_i)^{p^{m_i}}=\Big({\sum}_{h=1}^{p^{m_i}-1}\alpha(g_i,g_i^h),\ g_i^{p^{m_i}}\Big)=(p^{k_i},1).$$
Consequently, $\vartheta(u)^{-p^{m_2}}=\vartheta([u,v])^{p^{k_2}}$, and $\vartheta(v)^{-p^{m_1}}=\vartheta([u,v])^{p^{k_1}}$.

Thus, $\vartheta$ is well-defined.

Obviously, $|(m_2,m_1,d;k_2,k_1)|=p^{d+m_1+m_2}=|H_{k_1,k_2}|$, and $\vartheta$ is surjective. Therefore, $\vartheta$ is an isomorphism.
\end{proof}

In this way, we recover \cite[Theorem 1.1]{AMM12} in the case $p>2$.

\subsection{Orders of automorphism groups}

Denote $\mathcal{A}_{(p^{k_1},p^{k_2}),1}$, $\mathcal{L}_{(p^{k_1},p^{k_2}),1}$ respectively as $\mathcal{A}(k_1,k_2)$, $\mathcal{L}(k_1,k_2)$.

By Theorem \ref{thm:auto-group}, there is a short exact sequence
$$1\to \mathbb{Z}_{p^d}^2\to{\rm Aut}(H_{k_1,k_2})\to \mathcal{L}(k_1,k_2)\to 1.$$
Recall that $\mathbf{x}=(x_{ij})\in\mathcal{X}_G$ satisfies $p^{m_2-m_1}|x_{21}$ and $p\nmid\det(\mathbf{x})$.
Remembering that $\varphi_{\mathbf{x}}=\varphi_{\mathbf{x}'}$ if and only if $\mathbf{x}_{i,j}\equiv\mathbf{x}'_{i,j}\pmod{p^{m_i}}$ for all $i,j$, we may just assume $x_{11},x_{12}\in\{0,1,\ldots,p^{m_1}-1\}$ and $x_{21},x_{22}\in\{0,1,\ldots,p^{m_2}-1\}$.

By (\ref{eq:def-A}), $\mathbf{x}\in\mathcal{A}(k_1,k_2)$ if and only if $\mathbf{a}\mathbf{x}^\star=\det(\mathbf{x})\cdot\mathbf{a}$, which reads
\begin{align}
p^{k_2}\frac{x_{21}}{p^{m_2-m_1}}&=p^{k_1}(x_{11}x_{22}-x_{12}x_{21}-x_{11}),   \label{eq:auto-2-1}  \\
p^{k_1+m_2-m_1}x_{12}&=p^{k_2}(x_{11}x_{22}-x_{12}x_{21}-x_{22}),     \label{eq:auto-2-2}
\end{align}
as equations in $\mathbb{Z}_{p^d}$.

When $m_1=m_2=m$, $k_1=d$, $k_2=k$, (\ref{eq:auto-2-1}), (\ref{eq:auto-2-2}) are equivalent to
$$\|x_{21}\|\ge d-k,  \qquad  \|x_{11}-1\|\ge d-k.$$
\begin{itemize}
  \item If $k=d$, then $\mathcal{L}(k_1,k_2)={\rm GL}(2,\mathbb{Z}_{p^m})$. The reduction map $\mathbb{Z}_{p^m}\twoheadrightarrow\mathbb{Z}_{p}$
        defines an epimorphism
        $f:{\rm GL}(2,\mathbb{Z}_{p^m})\twoheadrightarrow{\rm GL}(2,\mathbb{Z}_{p})$, with $|\ker(f)|=p^{4(m-1)}$. As is well-known,
        $|{\rm GL}(2,\mathbb{Z}_p)|=(p^2-1)(p^2-p)$, so
        $$|\mathcal{L}(k_1,k_2)|=p^{4m-3}(p-1)(p^2-1).$$
  \item If $k<d$, then we can write $x_{11}=1+p^{d-k}u$, $x_{21}=p^{d-k}v$ for $u,v\in\{0,1,\ldots,p^{m-d+k}-1\}$; the invertibility of
        $\mathbf{x}$ is equivalent to $p\nmid x_{22}$, and there is no constraint on $x_{12}$. Hence
        $$|\mathcal{L}(k_1,k_2)|=p^{2m-2d+2k}\cdot(p^m-p^{m-1})\cdot p^m=p^{4m-2d+2k-1}(p-1).$$
\end{itemize}

Now suppose $m_2>m_1$. Then $p\nmid\det(\mathbf{x})$ is equivalent to $p\nmid x_{11}x_{22}$.
\begin{enumerate}
  \item When $k_1=k_2=d$, (\ref{eq:auto-2-1}), (\ref{eq:auto-2-2}) are trivial, i.e. there is no more constraint on the $x_{ij}$'s,
        so the numbers of choices for $x_{11},x_{12},x_{21},x_{22}$ are respectively $p^{m_1-1}(p-1)$, $p^{m_1}$, $p^{m_1}$, $p^{m_2-1}(p-1)$, implying
        $$|\mathcal{L}(k_1,k_2)|=p^{3m_1+m_2-2}(p-1)^2.$$
  \item When $0\le k_2=k<d=k_1$, (\ref{eq:auto-2-1}) is equivalent to
        $$\|x_{21}\|\ge m_2-m_1+d-k,$$
        and then (\ref{eq:auto-2-2}) becomes $\|x_{11}-1\|\ge d-k$.
        Hence the numbers of choices for $x_{11},x_{12},x_{21},x_{22}$ are respectively $p^{m_1-d+k}$, $p^{m_1}$, $p^{m_1-d+k}$, $p^{m_2-1}(p-1)$, implying
        $$|\mathcal{L}(k_1,k_2)|=p^{3m_1+m_2-2d+2k-1}(p-1).$$
  \item When $0\le k_1=k<d=k_2$, (\ref{eq:auto-2-2}) is equivalent to
        $$\|x_{12}\|\ge d-m_2+m_1-k,$$
        and then (\ref{eq:auto-2-1}) becomes $\|x_{22}-1\|\ge d-k$.
        The numbers of choices for $x_{11},x_{12},x_{21},x_{22}$ are respectively $p^{m_1-1}(p-1)$, $p^{m_2-d+k}$, $p^{m_1}$, $p^{m_2-d+k}$, implying
        $$|\mathcal{L}(k_1,k_2)|=p^{2m_1+2m_2-2d+2k-1}(p-1).$$
  \item When $0\le k_1<k_2<d\le k_1+m_2-m_1$, (\ref{eq:auto-2-2}) is equivalent to
        $$\|x_{11}-1\|\ge d-k_2.$$
        Clearly, $x_{21}\equiv p^{m_2-m_1}x_{11}z\pmod{p^d}$ for a unique $0\le z<p^{m_1}$, then (\ref{eq:auto-2-1}) is equivalent to
        $$\|x_{22}-1-p^{k_2-k_1}z\|\ge d-k_1.$$
        The numbers of choices for $x_{11},x_{12},x_{21},x_{22}$ are respectively $p^{m_1-d+k_2}$, $p^{m_1}$, $p^{m_1}$, $p^{m_2-d+k_1}$, implying
        $$|\mathcal{L}(k_1,k_2)|=p^{3m_1+m_2-2d+k_1+k_2}.$$
  \item When $0\le k_1<k_2<k_1+m_2-m_1<d$, there uniquely exist $z,z'\in\{0,1,\ldots,p^{m_1}-1\}$ such that
        $x_{21}\equiv p^{m_2-m_1}x_{11}z\pmod{p^{m_2}}$ and $x_{12}\equiv x_{22}z'\pmod{p^{m_1}}$,
        so (\ref{eq:auto-2-1}), (\ref{eq:auto-2-2}) are equivalent to
        \begin{align*}
        \|x_{11}-(1-p^{m_2-m_1}zz')^{-1}(1+p^{k_1+m_2-m_1-k_2}z')\|\ge d-k_2,  \\
        \|x_{22}-(1-p^{m_2-m_1}zz')^{-1}(1+p^{k_2-k_1}z)\|\ge d-k_1;
        \end{align*}
        the inverses are taken in $\mathbb{Z}_{p^d}$.
        Similarly as the previous case,
        $$|\mathcal{L}(k_1,k_2)|=p^{3m_1+m_2-2d+k_1+k_2}.$$
\end{enumerate}
Unify the results in Case 4 and Case 5 as: $|\mathcal{L}(k_1,k_2)|=p^{3m_1+m_2-2d+k_1+k_2}$ when $k_1<k_2<\min\{d,k_1+m_2-m_1\}$.

Multiplying by $|\mathbb{Z}_{p^d}^2|=p^{2d}$ yields the formula for the order of ${\rm Aut}(H_{k_1,k_2})$:
\begin{thm}\label{thm:order}
When $m_1=m_2=m$,
$$|{\rm Aut}(H_{d,k})|=\begin{cases} p^{4m+2d-3}(p-1)(p^2-1), &k=d \\ p^{4m+2k-1}(p-1),&k<d\end{cases};$$
when $m_1<m_2$,
$$|{\rm Aut}(H_{k_1,k_2})|=\begin{cases}
p^{3m_1+m_2+2d-2}(p-1)^2, &k_1=k_2=d \\
p^{3m_1+m_2+2k-1}(p-1),&k_2=k<d=k_1 \\
p^{2m_1+2m_2+2k-1}(p-1),&k_1=k<d=k_2 \\
p^{3m_1+m_2+k_1+k_2},&k_1<k_2<\min\{d,k_1+m_2-m_1\}
\end{cases}.$$
\end{thm}

\section{A family of Miller $p$-groups of minimal order}

Let $p>2$. In 1994, Morigi \cite{Mo94} showed that the order of a Miller $p$-group is at least $p^7$, and constructed the first Miller group of order $p^7$. Checking through the literature, we do not find the second one. In this section, we construct a family of new Miller $p$-groups of order $p^7$.

Let $r=3$, $n=4$, $d=m=1$, i.e. $A=\mathbb{Z}_p^3$, $G=\mathbb{Z}_p^4$. Let
$$\mathbf{a}=\left(\begin{array}{cccc} 1 & 0 & 0 & 0  \\ 0 & 1 & 0 & 0 \\ 0 & 0 & 1 & 0 \end{array}\right),  \qquad
\mathbf{u}=\left(\begin{array}{ccc} 1 & 1 & 0 \\ 0 & 1 & 1 \\ 0 & 0 & 1 \end{array}\right).$$

Enumerate the elements of $\Lambda$ as
\begin{align}
(2,3),\ \ \ (1,3), \ \ \ (1,2), \ \ \ (1,4), \ \ \ (2,4), \ \ \ (3,4).   \label{eq:enumeration}
\end{align}

Call $\vec{\mu}=(\mu_1,\mu_2,\mu_3)\in(\mathbb{Z}_p\setminus\{0\})^3$ {\it admissible} if $\mu_1,\mu_2,\mu_3$ are pairwise distinct.
For each admissible $\vec{\mu}$, put
$$\mathbf{c}(\vec{\mu})=(\mathbf{u},\mathbf{d}(\vec{\mu}))
=\left(\begin{array}{cccccc} 1 & 1 & 0 & \mu_1 & 0 & 0  \\ 0 & 1 & 1 & 0 & \mu_2 & 0 \\ 0 & 0 & 1 & 0 & 0 & \mu_3 \end{array}\right),$$
where $\mathbf{d}(\vec{\mu})$ denotes the diagonal matrix with entries $\mu_1,\mu_2,\mu_3$ on the main diagonal.
Denote $E(\alpha_{\mathbf{a},\mathbf{c}(\vec{\mu})})$ by $M_{\mu_1,\mu_2,\mu_3}$.

\begin{thm}\label{thm:Miller}
Each $M_{\mu_1,\mu_2,\mu_3}$ is a Miller group, and $M_{\mu_1,\mu_2,\mu_3}\cong M_{\nu_1,\nu_2,\nu_3}$ if and only if
$\mu_1/\nu_1=\mu_2/\nu_2=\mu_3/\nu_3$.
\end{thm}

\begin{proof}
Given admissible $\vec{\mu}=(\mu_1,\mu_2,\mu_3)$, $\vec{\nu}=(\nu_1,\nu_2,\nu_3)$,
consider
$$\mathbf{y}\mathbf{a}=\mathbf{a}\mathbf{x}, \qquad   \mathbf{y}\mathbf{c}(\vec{\mu})=\mathbf{c}(\vec{\nu})\mathbf{x}^\diamond,$$
for $\mathbf{y}=(y_{ij})\in{\rm GL}(3,\mathbb{Z}_p)$, $\mathbf{x}=(x_{ij})\in{\rm GL}(4,\mathbb{Z}_p)$.
We are going to show that these equations have no solution unless $\vec{\mu}$ is a multiple of $\vec{\nu}$, in which case there is a unique solution such that $\mathbf{y}$ is the identity matrix.

According to Convention \ref{conv:enumeration} and the definition (\ref{eq:x-diamond}), $(\mathbf{x}^\diamond)_{s,t}=x_{ui}x_{vj}-x_{uj}x_{vi}$ if the $s$-th and $t$-th elements in (\ref{eq:enumeration}) are respectively $(u,v)$, $(i,j)$.

Clearly, $\mathbf{y}\mathbf{a}=\mathbf{a}\mathbf{x}$ is equivalent to
$x_{14}=x_{24}=x_{34}=0$ and $x_{ij}=y_{ij}$ for $1\le i,j\le 3$.
When these hold, in block form
$\mathbf{x}^\diamond=\left(\begin{array}{cc} \mathbf{y}^\diamond & \mathbf{0} \\ \mathbf{z} & x_{44}\mathbf{y} \end{array}\right)$,
with
\begin{align*}
\mathbf{y}^\diamond&=\left(\begin{array}{ccc}
y_{22}y_{33}-y_{23}y_{32} & y_{21}y_{33}-y_{23}y_{31} & y_{21}y_{32}-y_{22}y_{31} \\
y_{12}y_{33}-y_{13}y_{32} & y_{11}y_{33}-y_{13}y_{31} & y_{11}y_{32}-y_{12}y_{31} \\
y_{12}y_{23}-y_{13}y_{22} & y_{11}y_{23}-y_{13}y_{21} & y_{11}y_{22}-y_{12}y_{21}
\end{array}\right),  \\
\mathbf{z}&=\left(\begin{array}{ccc}
y_{12}x_{43}-y_{13}x_{42} & y_{11}x_{43}-y_{13}x_{41} & y_{11}x_{42}-y_{12}x_{41} \\
y_{22}x_{43}-y_{23}x_{42} & y_{21}x_{43}-y_{23}x_{41} & y_{21}x_{42}-y_{22}x_{41} \\
y_{32}x_{43}-y_{33}x_{42} & y_{31}x_{43}-y_{33}x_{41} & y_{31}x_{42}-y_{32}x_{41}
\end{array}\right).
\end{align*}
Then $\mathbf{y}\mathbf{c}(\vec{\mu})=\mathbf{c}(\vec{\nu})\mathbf{x}^\diamond$ becomes
$$\mathbf{y}\mathbf{u}=\mathbf{u}\mathbf{y}^\diamond+\mathbf{d}(\vec{\nu})\mathbf{z}, \qquad
\mathbf{y}\mathbf{d}(\vec{\mu})=x_{44}\mathbf{d}(\vec{\nu})\mathbf{y}.$$

Let $\eta_i=x_{44}\nu_i$. It follows from $\mathbf{y}\mathbf{d}(\vec{\mu})=x_{44}\mathbf{d}(\vec{\nu})\mathbf{y}$ that
$(\eta_i-\mu_j)y_{ij}=0$ for all $i,j$. So $y_{ij}\ne0$ only if $\eta_i=\mu_j$.
Since $\det(\mathbf{y})\ne 0$ and $\vec{\mu},\vec{\nu}$ are admissible, there is a permutation $\tau$ on $\{1,2,3\}$ such that $\eta_i=\mu_{\tau(i)}$, $1\le i\le 3$.
\begin{enumerate}
  \item If $\eta_1=\mu_2$, $\eta_2=\mu_1$, $\eta_3=\mu_3$, then $y_{12}y_{21}y_{33}\ne 0$ and all the other $y_{ij}$'s vanish.
        As is easy to see, $(\mathbf{y}\mathbf{u})_{3,1}=(\mathbf{u}\mathbf{y}^\diamond)_{3,1}=0$.
        It follows from
        $$0=(\mathbf{y}\mathbf{u}-\mathbf{u}\mathbf{y}^\diamond-\mathbf{d}(\vec{\nu})\mathbf{z})_{3,1}=\nu_3y_{33}x_{42}$$
        that $x_{42}=0$, but then $(\mathbf{y}\mathbf{u}-\mathbf{u}\mathbf{y}^\diamond-\mathbf{d}(\vec{\nu})\mathbf{z})_{2,3}=y_{12}y_{21}\ne0$.
  \item If $\eta_1=\mu_2$, $\eta_2=\mu_3$, $\eta_3=\mu_1$, then $y_{12}y_{23}y_{31}\ne 0$ and all the other $y_{ij}$'s vanish.
        As is easy to see, $(\mathbf{y}\mathbf{u})_{1,1}=(\mathbf{u}\mathbf{y}^\diamond)_{1,1}=0$.
        It follows from
        $$0=(\mathbf{y}\mathbf{u}-\mathbf{u}\mathbf{y}^\diamond-\mathbf{d}(\vec{\nu})\mathbf{z})_{1,1}=-\nu_1y_{12}x_{43}$$
        that $x_{43}=0$, but then $(\mathbf{y}\mathbf{u}-\mathbf{u}\mathbf{y}^\diamond-\mathbf{d}(\vec{\nu})\mathbf{z})_{3,2}=y_{31}\ne0$.
  \item If $\eta_1=\mu_3$, $\eta_2=\mu_1$, $\eta_3=\mu_2$, then $y_{13}y_{21}y_{32}\ne 0$ and all the other $y_{ij}$'s vanish.
        As is easy to see, $(\mathbf{y}\mathbf{u})_{1,2}=(\mathbf{u}\mathbf{y}^\diamond)_{1,2}=0$.
        It follows from
        $$0=(\mathbf{y}\mathbf{u}-\mathbf{u}\mathbf{y}^\diamond-\mathbf{d}(\vec{\nu})\mathbf{z})_{1,2}=\nu_1y_{13}x_{41}$$
        that $x_{41}=0$, but then $(\mathbf{y}\mathbf{u}-\mathbf{u}\mathbf{y}^\diamond-\mathbf{d}(\vec{\nu})\mathbf{z})_{3,3}=y_{32}\ne0$.
  \item If $\eta_1=\mu_3$, $\eta_2=\mu_2$, $\eta_3=\mu_1$, then $y_{13}y_{22}y_{31}\ne 0$ and all the other $y_{ij}$'s vanish.
        As is easy to see,
        \begin{align*}
        \mathbf{y}^\diamond&=\left(\begin{array}{ccc} 0 & 0 & -y_{22}y_{31} \\
        0  & -y_{13}y_{31} & 0 \\  -y_{13}y_{22} & 0 & 0 \end{array}\right), \\
        \mathbf{z}&=\left(\begin{array}{ccc} -y_{13}x_{42} & -y_{13}x_{41} & 0 \\
        y_{22}x_{43}  & 0 & -y_{22}x_{41}  \\  0 & y_{31}x_{43} & y_{31}x_{42} \end{array}\right).
        \end{align*}
        So $(\mathbf{y}\mathbf{u})_{1,1}=(\mathbf{u}\mathbf{y}^\diamond)_{1,1}=0$.
        It follows from
        $$0=(\mathbf{y}\mathbf{u}-\mathbf{u}\mathbf{y}^\diamond-\mathbf{d}(\vec{\nu})\mathbf{z})_{1,1}=\nu_1y_{13}x_{42}$$
        that $x_{42}=0$. Then $\mathbf{y}\mathbf{u}-\mathbf{u}\mathbf{y}^\diamond-\mathbf{d}(\vec{\nu})\mathbf{z}=0$ reads
        $$\left(\begin{array}{ccc} 0 & y_{13}(y_{31}+\nu_1x_{41}) & y_{13}+y_{22}y_{31} \\
        y_{22}(y_{13}-\nu_2x_{43})  & y_{22}+y_{13}y_{31} & y_{22}(1+\nu_2x_{41}) \\
        y_{31}+y_{13}y_{22} & y_{31}(1-\nu_3x_{43}) & 0 \end{array}\right)=0.$$
        From the vanishing of the $(1,2)$- and $(2,3)$-entries we see $y_{31}=\nu_1/\nu_2$.
        From the vanishing of the $(2,1)$- and $(3,2)$-entries we see $y_{13}=\nu_2/\nu_3$.
        Then the vanishing of the anti-diagonal entries leads to $\nu_1^2=\nu_2^2=\nu_3^2$, which contradicts the assumption that $\nu_1,\nu_2,\nu_3$ are pairwise distinct.
  \item If $\eta_1=\mu_1$, $\eta_2=\mu_3$, $\eta_3=\mu_2$, then $y_{11}y_{23}y_{32}\ne 0$ and all the other $y_{ij}$'s vanish.
        As is easy to see, $(\mathbf{y}\mathbf{u})_{2,1}=(\mathbf{u}\mathbf{y}^\diamond)_{2,1}=0$.
        It follows from
        $$0=(\mathbf{y}\mathbf{u}-\mathbf{u}\mathbf{y}^\diamond-\mathbf{d}(\vec{\nu})\mathbf{z})_{2,1}=\nu_2y_{23}x_{42}$$
        that $x_{42}=0$, but then $(\mathbf{y}\mathbf{u}-\mathbf{u}\mathbf{y}^\diamond-\mathbf{d}(\vec{\nu})\mathbf{z})_{1,3}=-y_{11}y_{32}\ne0$.
\end{enumerate}
Thus, actually $\mu_i=\eta_i=x_{44}\nu_i$ for all $i$, implying that $\mathbf{y}$ is diagonal.

Suppose $\mathbf{y}=\mathbf{d}(\lambda_1,\lambda_2,\lambda_3)$, the diagonal matrix with the $(i,i)$-entry $\lambda_i$.
Then
$\mathbf{y}^\diamond=\mathbf{d}(\lambda_2\lambda_3,\lambda_1\lambda_3,\lambda_1\lambda_2)$,
and $\mathbf{y}\mathbf{u}-\mathbf{u}\mathbf{y}^\diamond-\mathbf{d}(\vec{\nu})\mathbf{z}=0$ reads
$$\left(\begin{array}{ccc} \lambda_1-\lambda_2\lambda_3 & \lambda_1(1-\lambda_3-\nu_1x_{43}) & -\nu_1\lambda_1x_{42} \\
-\nu_2\lambda_2x_{43} & \lambda_2-\lambda_1\lambda_3 & \lambda_2(1-\lambda_1+\nu_2x_{41})  \\
\nu_3\lambda_3x_{42} & \nu_3\lambda_3x_{41} & \lambda_3-\lambda_1\lambda_2 \end{array}\right)=0,$$
from which we can easily deduce
$$x_{41}=x_{42}=x_{43}=0,  \qquad  \lambda_1=\lambda_2=\lambda_3=1.$$

Thus, $M_{\mu_1,\mu_2,\mu_3}\cong M_{\nu_1,\nu_2,\nu_3}$ if and only if $\mu_1/\nu_1=\mu_2/\nu_2=\mu_3/\nu_3$.

In particular, when $(\mu_1,\mu_2,\mu_3)=(\nu_1,\nu_2,\nu_3)$, the only solution is $\mathbf{x},\mathbf{y}$ being the identity matrices, i.e. $\mathcal{L}_{\mathbf{a},\mathbf{c}}=1$. By Theorem \ref{thm:auto-group}, ${\rm Aut}(M_{\mu_1,\mu_2,\mu_3})\cong\hom(\mathbb{Z}_p^4,\mathbb{Z}_p^3)$ is abelian. Therefore, $M_{\mu_1,\mu_2,\mu_3}$ is a Miller group.
\end{proof}

As a conclusion, the $M_{\mu_1,\mu_2,1}$'s for distinct $\mu_1,\mu_2\in\mathbb{Z}_p\setminus\{0,1\}$ form a family of $(p-2)(p-3)$ pairwise non-isomorphic Miller groups of order $p^7$.

\begin{rmk}
\rm The group $G(1)$ given on \cite[Page 12]{Mo94} was defined as the group of class 2 generated by $a_1,a_2,b_1,b_2$, subject to the following relations:
\begin{align*}
a_1^p=a_2^p=1, \qquad [a_1,b_1]=[a_2,b_2]=[b_1,b_2]=1,  \\
b_1^p=[a_1,a_2][a_1,b_2][a_2,b_1], \qquad  b_2^p=[a_1,a_2][a_2,b_1].
\end{align*}
It is isomorphic to $E(\alpha_{\mathbf{a}_0,\mathbf{c}_0})$ with
$$\mathbf{a}_0=\left(\begin{array}{cccc} 0 & 0 & 1 & 1  \\ 0 & 0 & 1 & 1 \\ 0 & 0 & 1 & 0 \end{array}\right), \qquad
\mathbf{c}_0=\left(\begin{array}{cccccc} 1 & 0 & 0 & 0 & 0 & 0  \\ 0 & 0 & 1 & 0 & 0 & 0 \\ 0 & 0 & 0 & 1 & 0 & 0 \end{array}\right),$$
through
$a_1\mapsto(0,g_1)$, $a_2\mapsto(0,g_2)$, $b_1\mapsto(0,g_3)$, $b_2\mapsto(0,g_4).$

If $G(1)\not\cong M_{\mu_1,\mu_2,\mu_3}$ for some admissible $\vec{\mu}$, then $E(\alpha_{\mathbf{a}_0,\mathbf{c}_0})\cong E(\alpha_{\mathbf{a},\mathbf{c}(\vec{\mu})})$, so there exists $\mathbf{y}\in{\rm GL}(3,\mathbb{Z}_p)$, $\mathbf{x}\in{\rm GL}(4,\mathbb{Z}_p)$
with $\mathbf{y}\mathbf{a}_0=\mathbf{a}\mathbf{x}$, but this would imply $\mathbf{x}_{i,j}=0$ for $i=1,2,3$, $j=1,2$,
contradicting $\det(\mathbf{x})\ne 0$.

Thus, $G(1)\not\cong M_{\mu_1,\mu_2,\mu_3}$ for any admissible $\vec{\mu}$.
\end{rmk}

\bigskip

\noindent
{\bf Acknowledgement}

The author would like to thank the anonymous referee for giving many constructive suggestions and comments to improve the paper.

\section*{Declarations and statements}

{\bf Funding}: No funding was received for conducting this study.

\noindent
{\bf Conflict of interest}: The author has no conflict of interest to declare, 
and has no relevant financial or non-financial interests to disclose.

\noindent
{\bf Data availability statement}: Data sharing not applicable to this article as no data sets were generated or analyzed during the current study.

\medskip

\noindent
Haimiao Chen (orcid: 0000-0001-8194-1264)\ \ \ \ {\it chenhm@math.pku.edu.cn} \\
Department of Mathematics, Beijing Technology and Business University, \\
Liangxiang Higher Education Park, Fangshan District, Beijing, China.

\end{document}